# THE SHANNON INFORMATION OF FILTRATIONS AND THE ADDITIONAL LOGARITHMIC UTILITY OF INSIDERS

BY STEFAN ANKIRCHNER, STEFFEN DEREICH AND PETER IMKELLER

*Humboldt-Universität zu Berlin, Technische Universität Berlin
and Humboldt-Universität zu Berlin*

The background for the general mathematical link between utility and information theory investigated in this paper is a simple financial market model with two kinds of small traders: less informed traders and insiders, whose extra information is represented by an enlargement of the other agents' filtration. The expected logarithmic utility increment, that is, the difference of the insider's and the less informed trader's expected logarithmic utility is described in terms of the information drift, that is, the drift one has to eliminate in order to perceive the price dynamics as a martingale from the insider's perspective. On the one hand, we describe the information drift in a very general setting by natural quantities expressing the probabilistic better informed view of the world. This, on the other hand, allows us to identify the additional utility by entropy related quantities known from information theory. In particular, in a complete market in which the insider has some fixed additional information during the entire trading interval, its utility increment can be represented by the Shannon information of his extra knowledge. For general markets, and in some particular examples, we provide estimates of maximal utility by information inequalities.

**0. Introduction.** A simple mathematical model of two agents on a financial market taking their portfolio decisions on the basis of different information horizons has attracted much attention in recent years. Both agents are small, and unable to influence the price dynamics of the risky assets constituting the market. One agent just acts on the basis of the evolution of the market, the other one, the insider, possesses some additional knowledge at every instant of the continuous trading interval. This basic fact is modeled by associating two different filtrations with each agent, from which they make their portfolio decisions: the less informed agent, at time $t$, just has the $\sigma$-field $\mathcal{F}_t$, corresponding to the natural evolution of the market up to this time, at his disposal for deciding about future investments, while the insider is able to make better decisions, taking his knowledge from a bigger $\sigma$-field $\mathcal{G}_t \supset \mathcal{F}_t$. We give a short selection from among the many papers dealing with this model, just indicating the most important mathematical techniques used for its investigation. Methods are focused on martingale and stochastic control theory, and techniques of enlargement of filtrations (see [21]), starting with the conceptual paper by Duffie and Huang [10], mostly in the *initial enlargement* setting, that is, the insider gets some fixed extra information at the beginning of the trading interval. The model is successively studied on stochastic bases with increasing complexity: that is, Karatzas and Pikovsky [23] on Wiener space, Grorud and Pontier [14] allow Poissonian noise and Biagini and Oksendal [6] employ anticipative calculus techniques. In the same setting, Amendinger, Becherer and Schweizer [1] calculate the value of insider information from the perspective of specific utilities. Baudoin [5] introduces the concept of weak additional information consisting in the knowledge of the law of some random element. Campi [7] considers hedging techniques for insiders in the incomplete market setting. It is clear that the expected utility the insider is able to gain from final wealth in this simple model will be bigger than the uninformed traders' utility, for every utility function. And, in fact, many of the quoted papers deal with the calculation of a better informed agent's additional utility.

In [3], in the setting of initial enlargements and logarithmic utility, a crucial and natural link between the additional expected logarithmic utility and information theoretic concepts was made. The insider's logarithmic utility advantage is identified with the Shannon entropy of the additional information. In the same setting, Gasbarra and Valkeila [15] extended this link by interpreting the logarithmic utility increment by the Kullback–Leibler information of the insider's additional knowledge from the perspective of Bayesian modeling. In the environment of this utility-information paradigm, the papers [2, 8, 17, 18, 19, 20] describe additional utility and treat arbitrage questions and their interpretation in information theoretic terms in increasingly complex models of the same base structure, including some simple examples of progressive enlargements. It is clear that utility concepts different from the logarithmic one correspond on the information theoretic side to the generalized entropy concepts of $f$-divergences.

In this paper we shall continue the investigation of mathematical questions related to the link between utility and information theory in the most general setting of enlargements of filtrations: besides assuming eventually that the base space be standard, to ensure the existence of regular conditional probabilities, we shall let the filtration of the better informed agent just contain the one of the natural evolution of knowledge. To concentrate on one kind of entropy in this general setting, we shall consider logarithmic utility throughout. In this framework, Ankirchner and Imkeller [2] calculate the maximal expected utility of an agent from the intrinsic point of view of his (general) filtration, and relate the finiteness of expected utility via the (NFLVR) condition to the characterization of semimartingales by the theorem of Dellacherie–Meyer–Mokobodski. The compensator in the Doob–Meyer decomposition of underlying asset price processes with respect to the agent's filtration is determined by the *information drift process*. In this paper we shall give a general analysis of the nature of this process, and relate it to measuring the difference of the information residing in the two filtrations, independently of the particular price dynamics. The basic observation we start with in Section 2 identifies the information drift process with Radon–Nikodym densities of the stochastic kernel in an integral representation of the conditional probability process and the conditional probability process itself. This observation allows for an identification of the additional utility by the information difference of the two filtrations in terms of Shannon entropy notions in Section 5, again independent of particular price dynamics of the financial market.

The paper is organized as follows. In the preparatory Section 1 we recall the main results about the connection between finite utility filtrations, properties of the price dynamics from the perspective of different agents and properties of the information drift from [2]. In Section 2 (Theorems 2.6 and 2.10) properties of the conditional probability processes with respect to the agents' filtrations and the information drift process are investigated in depth, and lead to the identification of the information drift by subjective conditional probability quantities. The description of the additional utility in terms of entropy notions is more easily obtained, if the additional information in the bigger filtration comes in discrete bits along a sequence of partitions of the trading interval, leading to stepwise "initial enlargements" which ultimately converge to the big filtration as the mesh of the partitions shrinks to 0. This is done in Section 5 (Theorem 5.8), after being prepared in Sections 3 and 4 by a general investigation of the convergence properties of information drifts going along with the convergence of such discretized enlargements to the big filtration. In the final Section 6 general facts known from Shannon information theory (see [16]) are applied to estimate the expected maximal logarithmic utility of a better informed agent via the identification theorem of Section 5, in several particular cases. Entropy maximizing properties of Gaussian random variables play an important role.

**1. Preliminaries.** In this preparatory section we define the financial market model and recall some basic facts about expected utility maximization. Our favorite utility function will be the logarithmic one, for which we will then compare the maximal expected utilities of agents on the market who act on the background of asymmetric information. Recalling a result from [2], we will describe the utility increment of a better informed agent by the respective *information drift* of the agents' filtrations.

Let $(\Omega, \mathcal{F}, P)$ be a probability space with a filtration $(\mathcal{F}_t)_{0 \leq t \leq T}$, where $T > 0$ is a fixed time horizon. We consider a financial market with one nonrisky asset of interest rate normalized to 0, and one risky asset with price $S_t$ at time $t \in [0, T]$. We assume that $S$ is a continuous $(\mathcal{F}_t)$-semimartingale with values in $\mathbb{R}$ and write $\mathcal{A}$ for the set of all $S$-integrable and $(\mathcal{F}_t)$-predictable processes such that $\theta_0 = 0$. If $\theta \in \mathcal{A}$, then we denote by $(\theta \cdot S)$ the usual stochastic integral process. For all $x > 0$, we interpret

$$x + (\theta \cdot S)_t, \qquad 0 \leq t \leq T,$$

as the wealth process of a trader possessing an initial wealth $x$ and choosing the investment strategy $\theta$ on the basis of his knowledge horizon corresponding to the filtration $(\mathcal{F}_t)$. Throughout this paper we will suppose the preferences of the agents to be described by the logarithmic utility function. Furthermore, we suppose that the traders' total wealth has always to be strictly positive, that is, for all $t \in [0, T]$,

(1) $$x + (\theta \cdot S)_t > 0 \qquad \text{a.s.}$$

Strategies $\theta$ satisfying equation (1) will be called *x-superadmissible*. The agents want to maximize their expected logarithmic utility from their wealth at time $T$. So we are interested in the exact value of

$$u(x) = \sup\{E \log(x + (\theta \cdot S)_T) : \theta \in \mathcal{A} \ x\text{-superadmissible}\}.$$

Sometimes we will write $u_{\mathcal{F}}(x)$ in order to stress the underlying filtration. The expected logarithmic utility of the agent can be calculated easily, if one has a semimartingale decomposition of the form

(2) $$S_t = M_t + \int_0^t \eta_s \, d\langle M, M \rangle_s,$$

where $\eta$ is a predictable process. Such a decomposition is given for a large class of semimartingales. For example, if $S$ satisfies the property (NFLVR), then it may be decomposed as in (2) (see [12]). As is shown in a forthcoming Ph.D. thesis [4], finiteness of $u(x)$ implies already such a decomposition to exist. Hence, a decomposition as in (2) may be given even in cases where arbitrage exists. We state Theorem 2.9 of [2].

PROPOSITION 1.1. *Suppose $S$ can be decomposed into $S = M + \eta \cdot \langle M, M \rangle$. Then for any $x > 0$, the following equation holds:*

(3) $$u(x) = \log x + \tfrac{1}{2} E \int_0^T \eta_s^2 \, d\langle M, M \rangle_s.$$

This proposition motivates the following definition.

DEFINITION 1.2. A filtration $(\mathcal{G}_t)$ is called a *finite utility filtration* for $S$, if $S$ is a $(\mathcal{G}_t)$-semimartingale with decomposition $dS = dM + \zeta \cdot d\langle M, M \rangle$, where $\zeta$ is $(\mathcal{G}_t)$-predictable and belongs to $L^2(M)$, that is, $E \int_0^T \zeta^2 \, d\langle M, M \rangle < \infty$. We write

$$\mathbb{F} = \{(\mathcal{H}_t) \supset (\mathcal{F}_t) | (\mathcal{H}_t) \text{ is a finite utility filtration for } S\}.$$

We now compare two traders who take their portfolio decisions not on the basis of the same filtration, but on the basis of different information flows represented by the filtrations $(\mathcal{G}_t)$ and $(\mathcal{H}_t)$, respectively. Suppose that both filtrations $(\mathcal{G}_t)$ and $(\mathcal{H}_t)$ are finite utility filtrations. We denote by

(4) $$S = M + \zeta \cdot \langle M, M \rangle$$

the semimartingale decomposition with respect to $(\mathcal{G}_t)$ and by

(5) $$S = N + \beta \cdot \langle N, N \rangle$$

the decomposition with respect to $(\mathcal{H}_t)$. Obviously,

$$\langle M, M \rangle = \langle S, S \rangle = \langle N, N \rangle$$

and, therefore, the utility difference is equal to

$$u_{\mathcal{H}}(x) - u_{\mathcal{G}}(x) = \tfrac{1}{2} E \int_0^T (\beta^2 - \zeta^2) \, d\langle M, M \rangle.$$

Furthermore, (4) and (5) imply

(6) $$M = N - (\zeta - \beta) \cdot \langle M, M \rangle \qquad \text{a.s.}$$

If $\mathcal{G}_t \subset \mathcal{H}_t$ for all $t \geq 0$, equation (6) can be interpreted as the semimartingale decomposition of $M$ with respect to $(\mathcal{H}_t)$. In this case one can show that the utility difference depends only on the process $\mu = \zeta - \beta$. We therefore use the following notion.

DEFINITION 1.3. Let $(\mathcal{G}_t)$ be a finite utility filtration and $S = M + \zeta \cdot \langle M, M \rangle$ the Doob–Meyer decomposition of $S$ with respect to $(\mathcal{G}_t)$. Suppose that $(\mathcal{H}_t)$ is a filtration such that $\mathcal{G}_t \subset \mathcal{H}_t$ for all $t \in [0, T]$. The $(\mathcal{H}_t)$-adapted measurable process $\mu$ satisfying

$$M - \int_0^{\cdot} \mu_t \, d\langle M, M \rangle_t \text{ is a } (\mathcal{H}_t)\text{-local martingale}$$

is called *information drift* (see [19]) of $(\mathcal{H}_t)$ with respect to $(\mathcal{G}_t)$.

The following proposition relates the information drift to the expected logarithmic utility increment.

PROPOSITION 1.4. *Let $(\mathcal{G}_t)$ and $(\mathcal{H}_t)$ be two finite utility filtrations such that $\mathcal{G}_t \subset \mathcal{H}_t$ for all $t \in [0,T]$. If $\mu$ is the information drift of $(\mathcal{H}_t)$ w.r.t. $(\mathcal{G}_t)$, then we have*

$$u_{\mathcal{H}}(x) - u_{\mathcal{G}}(x) = \tfrac{1}{2} E \int_0^T \mu^2 \, d\langle M, M \rangle.$$

PROOF. See Theorem 2.13 in [2]. □

So far we only required the information drift to be measurable and adapted. Due to the continuity of $S$, we have the following.

PROPOSITION 1.5. *The information drift, provided it exists, may be chosen to be predictable.*

PROOF. Suppose $\mu$ is a measurable and $(\mathcal{G}_t)$-adapted process such that

$$M - \int_0^{\cdot} \mu_t \, d\langle M, M \rangle_t$$

is a $(\mathcal{G}_t)$-local martingale. We denote by $^p\mu$ the predictable projection of $\mu$ with respect to $(\mathcal{G}_t)$. We will show that $M - {}^p\mu \cdot \langle M, M \rangle$ remains a $(\mathcal{G}_t)$-local martingale.

Let $\tau$ be a stopping time localizing $M$ such that $M^\tau$, the martingale $M$ stopped at $\tau$, is bounded. To simplify notation, we assume $M^\tau = M$. Let $0 \leq s < t$, $A \in \mathcal{G}_s$ and $\varepsilon > 0$. Then

$$\begin{aligned}
E(\mathbb{1}_A(M_t - M_{s+\varepsilon})) &= E\left( \mathbb{1}_A \int_{s+\varepsilon}^t \mu_r \, d\langle M, M \rangle_r \right) \\
&= E\left( \mathbb{1}_A E\left[ \int_{s+\varepsilon}^t \mu_r \, d\langle M, M \rangle_r \Big| \mathcal{G}_s \right] \right) \\
&= E\left( \mathbb{1}_A E\left[ \int_{s+\varepsilon}^t {}^p\mu_r \, d\langle M, M \rangle_r \Big| \mathcal{G}_s \right] \right) \\
&= E\left( \mathbb{1}_A \int_{s+\varepsilon}^t {}^p\mu_r \, d\langle M, M \rangle_r \right)
\end{aligned}$$

(see Theorem 57, Chapter VI in [11]). By dominated convergence, the left-hand side of this equation converges to $E(\mathbb{1}_A(M_t - M_s))$ as $\varepsilon \downarrow 0$. The right-hand side converges by similar arguments. Hence, we obtain

$$E(\mathbb{1}_A(M_t - M_s)) = E\left( \mathbb{1}_A \int_s^t {}^p\mu_r \langle M, M \rangle_r \right),$$

which means that $M - {}^p\mu \cdot \langle M, M \rangle$ is a $(\mathcal{G}_t)$-martingale. □

We close this section by recalling some basic properties of information drifts.

LEMMA 1.6. *Suppose the filtration $(\mathcal{F}_t)$ is a finite utility filtration with respect to which the Doob–Meyer decomposition of $S$ is given by $S = M + \eta \cdot \langle M, M \rangle$. Let $(\mathcal{H}_t)$ be a filtration satisfying $\mathcal{F}_t \subset \mathcal{H}_t$ for all $t \in [0, T]$ and suppose that $(\mathcal{H}_t)$ has an information drift $\mu$ with respect to $(\mathcal{F}_t)$. Then the following properties hold true:*

(i) *If $\mu$ belongs to $L^2(M)$, then the maximal expected utility $u_\mathcal{H}(x)$ is finite for all $x > 0$.*

(ii) *The set of finite utility filtrations $\mathbb{F}$ is equal to the set of all filtrations containing $(\mathcal{F}_t)$ and possessing an information drift $\lambda$ with respect to $(\mathcal{F}_t)$ such that $\lambda \in L^2(M)$.*

(iii) *If $(\mathcal{H}_t)$ is a finite utility filtration, then $\mu$ is orthogonal to $L^2_\mathcal{F}(M)$, the subspace of $(\mathcal{F}_t)$-predictable processes in $L^2(M)$.*

(iv) *If $(\mathcal{G}_t)$ is a filtration such that $\mathcal{F}_t \subset \mathcal{G}_t \subset \mathcal{H}_t$ for all $t \in [0, T]$, then there is also an information drift $\kappa$ of $(\mathcal{G}_t)$ with respect to $(\mathcal{F}_t)$. More precisely, $\kappa$ is equal to the $L^2(M)$-projection of $\mu$ onto the subspace of the $(\mathcal{G}_t)$-predictable processes.*

PROOF. Properties (i) and (ii) are obvious. For property (iii), let $S = N + \beta \cdot \langle N, N \rangle$ denote the Doob–Meyer decomposition of $S$ relative to $(\mathcal{H}_t)$, and let $\theta \in L^2_\mathcal{F}(M)$. Since $\theta$ is adapted to both $(\mathcal{F}_t)$ and $(\mathcal{H}_t)$, the integrals $(\theta \cdot M)$ and $(\theta \cdot N)$ are square integrable martingales with expectation zero. Therefore,

$$E \int_0^T \theta \mu \, d\langle M, M \rangle = E \left[ \int_0^T \theta \beta \, d\langle M, M \rangle - \int_0^T \theta \eta \, d\langle M, M \rangle \right]$$
$$= E \left[ \int_0^T \theta \, dM - \int_0^T \theta \, dN \right]$$
$$= 0.$$

Thus, $\mu$ is orthogonal to $L^2_\mathcal{F}(M)$. For property (iv), we refer again to [2]. □

**2. General enlargements.** Assume again that the price process $S$ is a semimartingale of the form

$$S = M + \eta \cdot \langle M, M \rangle,$$

with respect to a finite utility filtration $(\mathcal{F}_t)$. Moreover, let $(\mathcal{G}_t)$ be a filtration such that $\mathcal{F}_t \subset \mathcal{G}_t$, and let $\alpha$ be the information drift of $(\mathcal{G}_t)$ relative to $(\mathcal{F}_t)$. And, for simplicity of notation, suppose in this section that time horizon is infinite, that is, $T = \infty$. We shall aim at describing the relative information drift $\alpha$ by basic quantities related to the conditional probabilities of the larger $\sigma$-algebras $\mathcal{G}_t$ with respect to the smaller ones $\mathcal{F}_t, t \geq 0$.

Roughly, modulo some tedious technical details to be specified below, the relationship is as follows. Suppose for all $t \geq 0$ there is a regular conditional probability $P_t(\cdot, \cdot)$ of $\mathcal{F}$ given $\mathcal{F}_t$, which can be decomposed into a martingale component orthogonal to $M$, plus a component possessing a stochastic integral representation with respect to $M$ with a kernel function $k_t(\cdot, \cdot)$. Then we shall see that, provided $\alpha$ is square integrable with respect to $d\langle M, M\rangle \otimes P$, the kernel function at $t$ will be a signed measure in its set variable. Moreover, this measure is absolutely continuous with respect to the conditional probability, if restricted to $\mathcal{G}_t$, and $\alpha$ coincides with their Radon–Nikodym density.

We shall even be able to show that this relationship also makes sense in the reverse direction. Roughly, if absolute continuity of the stochastic integral kernel with respect to the conditional probabilities holds, and the Radon–Nikodym density is square integrable, the latter will turn out to provide an information drift $\alpha$ in a Doob–Meyer decomposition of $S$ in the larger filtration.

We shall finish the section with an illustration of this fundamental relationship by discussing some simple examples of particularly enlarged filtrations.

The discussion of the details of this fundamental relationship requires some care with the complexity of the underlying filtrations and state spaces. Of course, the need to work with conditional probabilities first of all confines us to spaces on which they exist. Let therefore $(\Omega, \mathcal{F}, P)$ be a standard Borel probability space (see [22]) with a filtration $(\mathcal{F}_t^0)_{t \geq 0}$ consisting of countably generated $\sigma$-algebras, and $M$ a $(\mathcal{F}_t^0)$-local martingale. We will also deal with the smallest right-continuous and completed filtration containing $(\mathcal{F}_t^0)$, which we denote by $(\mathcal{F}_t)$. We suppose that $\mathcal{F}_0$ is trivial and that every $(\mathcal{F}_t)$-local martingale has a continuous modification. Since $\mathcal{F}_t^0$ is a subfield of a standard Borel space, there exist regular conditional probabilities $P_t$ relative to the $\sigma$-algebras $\mathcal{F}_t^0$. Then for any set $A \in \mathcal{F}$, the process

$$(t, \omega) \mapsto P_t(\omega, A)$$

is an $(\mathcal{F}_t^0)$-martingale with a continuous modification (see, e.g., Theorem 4, Chapter VI in [11]). Note that the modification may not be adapted to $(\mathcal{F}_t^0)$, but only to $(\mathcal{F}_t)$. Furthermore, it is no problem to assume that the processes $P_t(\cdot, A)$ are modified in a way such that $P_t(\omega, \cdot)$ remains a measure on $\mathcal{F}$ for $P_M$-almost all $(\omega, t)$, where $P_M$ is a measure on $\Omega \times \mathbb{R}_+$ defined by $P_M(\Gamma) = E \int_0^\infty \mathbb{1}_\Gamma(\omega, t) \, d\langle M, M\rangle_t$, $\Gamma \in \mathcal{F} \otimes \mathcal{B}_+$.

It is known that each of these martingales may be uniquely written (see, e.g., [24], Chapter V)

$$(7) \qquad P_t(\cdot, A) = P(A) + \int_0^t k_s(\cdot, A) \, dM_s + L_t^A,$$

where $k(\cdot, A)$ is $(\mathcal{F}_t)$-predictable and $L^A$ satisfies $\langle L^A, M \rangle = 0$.

Now let $(\mathcal{G}_t^0)$ be another filtration on $(\Omega, \mathcal{F}, P)$ satisfying

$$\mathcal{F}_t^0 \subset \mathcal{G}_t^0$$

for all $0 \leq t \leq T$. We assume that each $\sigma$-field $\mathcal{G}_t^0$ is generated by a countable number of sets, and denote by $(\mathcal{G}_t)$ the smallest right-continuous and completed filtration containing $(\mathcal{G}_t^0)$. It is clear that each $\sigma$-field in the left-continuous filtration $(\mathcal{G}_{t-}^0)$ is also generated by a countable number of sets. We claim that the existence of an information drift of $(\mathcal{G}_t)$ relative to $(\mathcal{F}_t)$ for the process $M$ depends on whether the following condition is satisfied or not.

CONDITION 2.1. $k_t(\omega, \cdot)|_{\mathcal{G}_{t-}^0}$ is a signed measure and satisfies

$$k_t(\omega, \cdot)|_{\mathcal{G}_{t-}^0} \ll P_t(\omega, \cdot)|_{\mathcal{G}_{t-}^0}$$

for $P_M$-a.a $(\omega, t)$.

REMARK 2.2. Unfortunately, we have to distinguish between the filtrations $(\mathcal{F}_t^0)$, $(\mathcal{G}_t^0)$ and their extensions $(\mathcal{F}_t)$, $(\mathcal{G}_t)$. The reason is that the regular conditional probabilities considered exist only with respect to the smaller $\sigma$-fields. On the other hand, we use stochastic integration techniques which were developed only under the assumption that the underlying filtrations satisfy the usual conditions, and this necessitates working also with the larger $\sigma$-fields.

Let us next state some essential properties of the Radon–Nikodym density process existing according to our condition.

LEMMA 2.3. *Suppose Condition 2.1 is satisfied. Then there exists an $(\mathcal{F}_t \otimes \mathcal{G}_t)$-predictable process $\gamma$ such that, for $P_M$-a.a. $(\omega, t)$,*

$$\gamma_t(\omega, \omega') = \left.\frac{dk_t(\omega, \cdot)}{dP_t(\omega, \cdot)}\right|_{\mathcal{G}_{t-}^0}(\omega').$$

REMARK 2.4. Note that $\gamma_t(\omega, \cdot)$ is $\mathcal{G}_{t-}$-measurable. This is due to the fact that the predictable $\sigma$-algebra does not change by taking the left-continuous version of the underlying filtration.

PROOF OF LEMMA 2.3. Let $t_i^n = \frac{i}{2^n}$ for all $n \geq 0$ and $i \geq 0$. We denote by $\mathbb{T}$ the set of all $t_i^n$. It is possible to choose a family of finite partitions $(\mathcal{P}^{i,n})$ such that:

- for all $t \in \mathbb{T}$, we have $\mathcal{G}_{t-}^0 = \sigma(\mathcal{P}^{i,n} : i, n \geq 0 \text{ s.t. } t_i^n = t)$,

- $\mathcal{P}^{i,n} \subset \mathcal{P}^{i+1,n}$,
- if $i < j$, $n < m$ and $i2^{-n} = j2^{-m}$, then $\mathcal{P}^{i,n} \subset \mathcal{P}^{j,m}$.

We define, for all $n \geq 0$,

$$\gamma_t^n(\omega, \omega') = \sum_{i \geq 0} \sum_{A \in \mathcal{P}^{i,n}} \mathbb{1}_{]t_i^n, t_{i+1}^n]}(t) \mathbb{1}_A(\omega') \frac{k_t(\omega, A)}{P_t(\omega, A)}.$$

Note that $\frac{k_t(\omega, A)}{P_t(\omega, A)}$ is $(\mathcal{F}_t)$-predictable and $\mathbb{1}_{]t_i^n, t_{i+1}^n]}(t) \mathbb{1}_A(\omega')$ is $(\mathcal{G}_t)$-predictable. Hence, the product of both functions, defined as a function on $\Omega^2 \times \mathbb{R}_+$, is predictable with respect to $(\mathcal{F}_t \otimes \mathcal{G}_t)$. It follows that each $\gamma^n$ and, thus,

$$\gamma = \liminf_{n \to \infty} \gamma^n,$$

is $(\mathcal{F}_t \otimes \mathcal{G}_t)$-predictable.

Now fix $t \geq 0$. We claim that $k_t(\omega, \cdot) = \int_{\cdot} \gamma_t(\omega, \omega') P_t(\omega, d\omega')$ and, hence, that $\gamma_t(\omega, \cdot)$ is the density of $k_t(\omega, \cdot)$ with respect to $P_t(\omega, \cdot)$, $P_M$-a.s. For all $n \geq 0$, let $j = j(n)$ be the integer satisfying $t_j^n < t \leq t_{j+1}^n$ and denote by $\mathcal{Q}^n$ the corresponding partition $\mathcal{P}^{j,n}$. Observe that $(\mathcal{Q}^n)$ is an increasing sequence of partitions satisfying

$$\sigma(\mathcal{Q}^n : n \geq 0) = \mathcal{G}_{t-}^0$$

and, hence,

$$\gamma_t(\omega, \omega') = \liminf_n \gamma_t^n(\omega, \omega')$$

$$= \liminf_n \sum_{A \in \mathcal{Q}^n} \mathbb{1}_A(\omega') \frac{k_t(\omega, A)}{P_t(\omega, A)} = \left.\frac{dk_t(\omega, \cdot)}{dP_t(\omega, \cdot)}\right|_{\mathcal{G}_{t-}^0}. \qquad \Box$$

LEMMA 2.5. *If $(t, \omega, \omega') \mapsto \theta_t(\omega, \omega')$ is $(\mathcal{F}_t \otimes \mathcal{G}_t)$-predictable and bounded, then*

$$\int\!\!\int\!\!\int \theta_t(\omega, \omega') P_t(\omega, d\omega') \, d\langle M, M \rangle_t \, dP(\omega) = \int\!\!\int \theta_t(\omega, \omega) \, d\langle M, M \rangle_t \, dP(\omega).$$

PROOF. Let $0 \leq r < s$, $A \in \mathcal{F}_r$, $B \in \mathcal{G}_r$ and

$$\theta_t(\omega, \omega') = \mathbb{1}_{]r,s]}(t) \mathbb{1}_A(\omega) \mathbb{1}_B(\omega').$$

Then

$$\int\!\!\int\!\!\int \theta_t(\omega, \omega') P_t(\omega, d\omega') \, d\langle M, M \rangle_t \, dP(\omega)$$

$$= \int\!\!\int_r^s \mathbb{1}_A(\omega) P_t(\omega, B) \, d\langle M, M \rangle_t \, dP(\omega)$$

$$= \int\!\!\int_r^s \mathbb{1}_A(\omega) \mathbb{1}_B(\omega) \, d\langle M, M \rangle_t \, dP(\omega)$$

$$= \int\!\!\int \theta_t(\omega, \omega) \, d\langle M, M \rangle_t \, dP(\omega),$$

where the second equality holds due to results about optional projections (see Theorem 57, Chapter VI in [11]). By a monotone class argument, this can be extended to all bounded and $(\mathcal{F}_t \otimes \mathcal{G}_t)$-predictable processes. □

THEOREM 2.6. *Suppose Condition 2.1 is satisfied and $\gamma$ is as in Lemma 2.3. Then*

$$\alpha_t(\omega) = \gamma_t(\omega, \omega)$$

*is the information drift of $(\mathcal{G}_t)$ relative to $(\mathcal{F}_t)$.*

PROOF. Suppose $\tau$ to be a stopping time such that $M^\tau$ is a martingale. For $0 \leq s < t$ and $A \in \mathcal{G}_s^0$, we have to show

$$E[\mathbb{1}_A(M_t^\tau - M_s^\tau)] = E\left[\mathbb{1}_A \int_s^t \gamma_u(\omega,\omega)\,d\langle M,M\rangle_u^\tau\right].$$

For notational simplicity, write $M^\tau = M$ and observe

$$\begin{aligned}
E[\mathbb{1}_A(M_t - M_s)] &= E[P_t(\cdot,A)(M_t - M_s)] \\
&= E\left[(M_t - M_s)\int_0^t k_u(\cdot,A)\,dM_u\right] + E[(M_t - M_s)L_t^A] \\
&= E\left[\int_s^t k_u(\cdot,A)\,d\langle M,M\rangle_u\right] \\
&= E\left[\int_s^t \int_A \gamma_u(\omega,\omega')\,dP_u(\omega,d\omega')\,d\langle M,M\rangle_u\right] \\
&= E\left[\mathbb{1}_A(\omega)\int_s^t \gamma_u(\omega,\omega)\,d\langle M,M\rangle_u\right],
\end{aligned}$$

where we used Lemma 2.5 in the last equation. □

COROLLARY 2.7. *$(\mathcal{G}_t)$ is a finite utility filtration if and only if*

$$\int\!\!\int\!\!\int \gamma_t^2(\omega,\omega')P_t(\omega,d\omega')\,d\langle M,M\rangle_t\,dP(\omega) < \infty.$$

PROOF. This follows immediately from Lemma 2.5. □

We now look at the problem from the reverse direction. Starting with the assumption that $(\mathcal{G}_t)$ is a finite utility filtration, which amounts to $E\int_0^T \alpha^2\,d\langle M,M\rangle < \infty$, we show the validity of Condition 2.1.

In the sequel, $(\mathcal{G}_t)$ denotes a finite utility filtration and $\alpha$ its predictable information drift, that is,

(8) $$\widetilde{M} = M - \int_0^\cdot \alpha_t\,d\langle M,M\rangle_t$$

is a $(\mathcal{G}_t)$-local martingale. To prove the main results (Theorems 2.10 and 2.12), we need the following lemma.

LEMMA 2.8. *Let $0 \leq s < t$ and $\mathcal{P} = \{A_1, \ldots, A_n\}$ be a finite partition of $\Omega$ into $\mathcal{G}_s^0$-measurable sets. Then*

$$E \int_s^t \sum_{k=1}^n \left(\frac{k_u}{P_u}\right)^2 (\cdot, A_k) \mathbb{1}_{A_k} d\langle M, M \rangle_u \leq 4 E \left( \int_s^t \alpha_u^2 \, d\langle M, M \rangle_u \right) < \infty.$$

PROOF. Let $\mathcal{P} = \{A_1, \ldots, A_n\}$ be a finite $\mathcal{G}_s^0$-partition. An application of Itô's formula, in conjunction with (7) and (8), yields

$$\sum_{k=1}^n [\mathbb{1}_{A_k} \log P_s(\cdot, A_k) - \mathbb{1}_{A_k} \log P_t(\cdot, A_k)]$$

$$= \sum_{k=1}^n \left[ -\int_s^t \frac{1}{P_u(\cdot, A_k)} \mathbb{1}_{A_k} dP_u(\cdot, A_k) \right.$$

$$\left. + \frac{1}{2} \int_s^t \frac{1}{P_u(\cdot, A_k)^2} \mathbb{1}_{A_k} d\langle P(\cdot, A_k), P(\cdot, A_k) \rangle_u \right]$$

(9) $$= \sum_{k=1}^n \left[ -\int_s^t \frac{k_u}{P_u}(\cdot, A_k) \mathbb{1}_{A_k} d\widetilde{M}_u - \int_s^t \frac{k_u}{P_u}(\cdot, A_k) \mathbb{1}_{A_k} \alpha_u \, d\langle M, M \rangle_u \right.$$

$$- \int_s^t \frac{1}{P_u(\cdot, A_k)} \mathbb{1}_{A_k} dL_u^{A_k}$$

$$+ \frac{1}{2} \int_s^t \left(\frac{k_u}{P_u}\right)^2 (\cdot, A_k) \mathbb{1}_{A_k} d\langle M, M \rangle_u$$

$$\left. + \frac{1}{2} \int_s^t \frac{1}{P_u(\cdot, A_k)^2} \mathbb{1}_{A_k} d\langle L^{A_k}, L^{A_k} \rangle_u \right].$$

Note that $P_t(\cdot, A_k) \log P_t(\cdot, A_k)$ is a submartingale bounded from below for all $k$. Hence, the expectation of the left-hand side in the previous equation is at most 0.

A priori it is not clear whether

$$\sum_{k=1}^n \int_s^t \frac{k_u}{P_u}(\cdot, A_k) \mathbb{1}_{A_k} d\widetilde{M}_u$$

is integrable or not. Consider therefore, for all $\varepsilon > 0$, stopping times defined by

$$\tau_k^\varepsilon = \begin{cases} \infty, & \omega \notin A_k, \\ \inf\{t \geq s : P_t(\cdot, A_k) \leq \varepsilon\}, & \text{else,} \end{cases}$$

and
$$\tau^\varepsilon = \tau_1^\varepsilon \wedge \cdots \wedge \tau_n^\varepsilon.$$

Observe that $\tau^\varepsilon \to \infty$ as $\varepsilon \downarrow 0$ and that the stopped process
$$\sum_{k=1}^n \int_s^{t\wedge\tau^\varepsilon} \frac{k_u}{P_u}(\cdot, A_k) \mathbb{1}_{A_k}\, d\widetilde{M}_u$$

has expectation zero, since
$$E\left[\left(\int_s^{t\wedge\tau^\varepsilon} \sum_{k=1}^n \frac{k_u}{P_u}(\cdot, A_k)\mathbb{1}_{A_k}\, d\widetilde{M}_u\right)^2\right]$$
$$= E\left[\int_s^{t\wedge\tau^\varepsilon} \sum_{k=1}^n \left(\frac{k_u}{P_u}\right)^2(\cdot, A_k)\mathbb{1}_{A_k}\, d\langle M, M\rangle_u\right]$$
$$\leq \frac{1}{\varepsilon^2} E\left[\int_s^{t\wedge\tau^\varepsilon} \sum_{k=1}^n (k_u)^2(\cdot, A_k)\mathbb{1}_{A_k}\, d\langle M, M\rangle_u\right]$$
$$\leq \frac{1}{\varepsilon^2} E\left[\sum_{k=1}^n \int_s^t d\langle P(\cdot, A_k), P(\cdot, A_k)\rangle_u\right]$$
$$< \infty.$$

Similarly, one can show that the expectation of
$$\int_s^{t\wedge\tau^\varepsilon} \frac{1}{P_u(\cdot, A_k)} \mathbb{1}_{A_k}\, dL_u^{A_k}$$

vanishes. Consequently, we may deduce from (9) and the Kunita–Watanabe inequality
$$E \sum_{k=1}^n \frac{1}{2} \int_s^{t\wedge\tau^\varepsilon} \left(\frac{k_u}{P_u}\right)^2(\cdot, A_k)\mathbb{1}_{A_k}\, d\langle M, M\rangle_u$$
$$\leq E \sum_{k=1}^n \left[\int_s^{t\wedge\tau^\varepsilon} \frac{k_u}{P_u}(\cdot, A_k)\mathbb{1}_{A_k}\alpha_u\, d\langle M, M\rangle_u\right]$$
$$\leq E\left(\int_s^{t\wedge\tau^\varepsilon} \sum_{k=1}^n \left(\frac{k_u}{P_u}\right)^2(\cdot, A_k)\mathbb{1}_{A_k}\, d\langle M, M\rangle_u\right)^{1/2}$$
$$\times E\left(\int_s^{t\wedge\tau^\varepsilon} \alpha_u^2\, d\langle M, M\rangle_u\right)^{1/2},$$

which implies
$$E \int_s^{t\wedge\tau^\varepsilon} \sum_{k=1}^n \left(\frac{k_u}{P_u}\right)^2(\cdot, A_k)\mathbb{1}_{A_k}\, d\langle M, M\rangle_u \leq 4E\left(\int_s^{t\wedge\tau^\varepsilon} \alpha_u^2\, d\langle M, M\rangle_u\right).$$

Now the proof may be completed by a monotone convergence argument. □

Let $\mathbb{T}$ and $(\mathcal{P}^{i,n})_{i,n\geq 0}$ be a family of partitions as in the proof of Lemma 2.3. We define, for all $n \geq 0$,

$$Z^n_t(\omega,\omega') = \sum_{i\geq 0} \sum_{A\in\mathcal{P}^{i,n}} \mathbb{1}_{]t^n_i, t^n_{i+1}]}(t) \mathbb{1}_A(\omega') \frac{k_t(\omega,A)}{P_t(\omega,A)}.$$

Note that $Z^n$ is $(\mathcal{F}_t \otimes \mathcal{G}_t)$-predictable. We are now able to prove a converse statement to Theorem 2.6. Observe first the following:

LEMMA 2.9. *For $P_M$-almost all $(\omega,t) \in \Omega \times \mathbb{R}_+$, the discrete process $(Z^m_t(\omega,\cdot))_{m\geq 1}$ is an $L^2(P_t(\omega,\cdot))$-bounded martingale.*

PROOF. Every statement in the sequel is meant to hold for $P_M$-a.a. $(\omega,t) \in \Omega \times \mathbb{R}_+$.

Let $m \geq 0$, $l \geq 0$ and $j$ be the natural number such that $]t^{m+1}_l, t^{m+1}_{l+1}] \subset ]t^m_j, t^m_{j+1}]$. We start by proving that on $]t^{m+1}_l, t^{m+1}_{l+1}]$ we have

$$E^{P_t(\omega,\cdot)}[Z^{m+1}_t(\omega,\cdot)|\mathcal{P}^{j,m}] = Z^m_t(\omega,\cdot).$$

For this, let $B \in \mathcal{P}^{j,m}$ and $A_1, \ldots, A_k \in \mathcal{P}^{l,m+1}$ such that $A_1 \cup \cdots \cup A_k = B$. Note that

$$\begin{aligned}
E^{P_t(\omega,\cdot)}[\mathbb{1}_B(\cdot)Z^{m+1}_t(\omega,\cdot)] &= E^{P_t(\omega,\cdot)}\left[\sum_{i=1}^k \mathbb{1}_{A_i}(\cdot)\frac{k_t}{P_t}(\omega,A_i)\right] \\
&= \sum_{i=1}^k k_t(\omega,A_i) \\
&= k_t(\omega,B) \\
&= E^{P_t(\omega,\cdot)}[\mathbb{1}_B(\cdot)Z^m_t(\omega,\cdot)]
\end{aligned}$$

on $]t^{m+1}_l, t^{m+1}_{l+1}]$. Consequently, the process $(Z^m_t(\omega,\cdot))_{m\geq 1}$ is a martingale (with respect to a filtration depending on $t$). The martingale property implies that the sequence $\int (Z^n_t)^2(\omega,\omega')P_t(\omega,d\omega')$ is increasing and, hence, by monotone convergence,

$$\sup_n E \iint (Z^n_t)^2(\omega,\omega')P_u(\omega,d\omega')\,d\langle M,M\rangle_t$$
$$= E \int \sup_n \int (Z^n_t)^2(\omega,\omega')P_u(\omega,d\omega')\,d\langle M,M\rangle_t.$$

By Lemmas 2.8 and 2.5, we have

$$\sup_n E \iint (Z_u^n)^2(\omega, \omega') P_u(\omega, d\omega') \, d\langle M, M\rangle_u$$

$$= \sup_n E \int (Z_u^n)^2(\omega, \omega) \, d\langle M, M\rangle_u$$

$$= \sup_n E \sum_{i\geq 0} \int_{t_i^n}^{t_{i+1}^n} \sum_{A\in\mathcal{P}^{i,n}} \mathbb{1}_A(\omega) \left(\frac{k_t(\omega, A)}{P_t(\omega, A)}\right)^2 d\langle M, M\rangle_u$$

$$\leq 4E\left(\int \alpha_u^2 \, d\langle M, M\rangle_u\right) < \infty.$$

This shows that $(Z^n)_{n\geq 1}$ is an $L^2(P_t(\omega, \cdot))$-bounded martingale. $\square$

We now will show that $k$ can be chosen to be a signed measure. For this we identify $P_t(\omega, \cdot)$ with another measure on a countable generator of $\mathcal{G}_{t-}^0$. We then apply the result that two Banach space valued measures are equal, if they coincide on a generator stable for finite intersections.

THEOREM 2.10. *The kernel $k$ may be chosen such that*

$$\mathcal{G}_{t-}^0 \ni A \mapsto k_t(\omega, A) \in \mathbb{R}$$

*is a signed measure which is absolutely continuous with respect to $P_t(\omega, \cdot)|_{\mathcal{G}_{t-}^0}$, for $P_M$-a.a. $(\omega, t) \in \Omega \times [0, \infty)$. This means that Condition 2.1 is satisfied.*

PROOF. Lemma 2.9 implies that $(Z_t^m(\omega, \cdot))_{m\geq 1}$ is an $L^2(P_t(\omega, \cdot))$-bounded martingale and, hence, for a.a. fixed $(\omega, t)$, $(Z_t^m(\omega, \cdot))_{m\geq 1}$ possesses a limit $Z$. It can be chosen to be $(\mathcal{F}_t \otimes \mathcal{G}_t)$-predictable. Take, for example,

$$Z_t = \liminf_n (Z_t^n \vee 0) + \limsup_n (Z_t^n \wedge 0).$$

Now define a signed measure by

$$\tilde{k}_t(\omega, A) = \int \mathbb{1}_A(\omega') Z_t(\omega, \omega') \, dP_t(\omega, d\omega').$$

Observe that $\tilde{k}_t(\omega, \cdot)$ is absolutely continuous with respect to $P_t(\omega, \cdot)$ and that we have, for all $A \in \mathcal{P}^{j,m}$ with $j2^{-m} \leq t$,

$$\tilde{k}_t(\omega, A) = k_t(\omega, A)$$

for $P_M$-a.a. $(\omega, t) \in \Omega \times \mathbb{R}_+$. One may also interpret $\mathcal{G}_{t-}^0 \ni A \mapsto \tilde{k}_t(\omega, A)$, as an $L^2(M)$-valued measure. By applying the stochastic integral operator, we obtain an $L^2(\Omega)$-valued measure: $\mathcal{G}_{t-}^0 \ni A \mapsto \int_0^t \tilde{k}_s(\omega, A) \, dM_s$. Moreover,

$$(10) \qquad P_t(\omega, A) = P(A) + \int_0^t \tilde{k}_s(\omega, A) \, dM_s + L_t^A(\omega)$$

for all $A \in \bigcup_{j2^{-m} \leq t} \mathcal{P}^{j,m}$. Since the LHS and both expressions on the RHS are measures coinciding on a system which is stable for intersections, (10) holds for all $A \in \mathcal{G}_{t-}^0$. Hence, by choosing $k_t(\cdot, A) = \tilde{k}_t(\cdot, A)$ for all $A \in \mathcal{G}_{t-}^0$, the proof is complete. $\square$

REMARK 2.11. Since $k$ is determined up to $P_M$-null sets, we may assume that $k_t(\omega, \cdot)$ is absolutely continuous relative to $P_t(\omega, \cdot)$ everywhere.

We close this section with some examples showing how (well-known) information drifts can be derived explicitly, based on the formalism of Theorem 2.6. To this end, it is not always necessary to determine the signed measures $k_t(\omega, \cdot)$ on the whole $\sigma$-algebras $\mathcal{G}_t^0$, but only on some sub-$\sigma$-fields. This is the case, for example, if

$$\mathcal{G}_t^0 = \mathcal{F}_t^0 \vee \mathcal{H}_t^0, \qquad 0 \leq t \leq T,$$

where $(\mathcal{H}_t^0)$ is some countably generated filtration on $(\Omega, \mathcal{F})$.

Now suppose that $k_t(\omega, \cdot)$ is a signed measure on $(\mathcal{H}_{t-}^0)$ satisfying

$$k_t(\omega, \cdot)|_{\mathcal{H}_{t-}^0} \ll P_t(\omega, \cdot)|_{\mathcal{H}_{t-}^0}$$

for $P_M$-a.a $(\omega, t)$. Then we can show with the arguments of the proof of Lemma 2.3 that there is an $(\mathcal{F}_t \otimes \mathcal{H}_t)$-predictable process $\beta$ such that $P_M$-a.e.

$$\beta_t(\omega, \omega') = \left. \frac{dk_t(\omega, \cdot)}{dP_t(\omega, \cdot)} \right|_{\mathcal{H}_{t-}^0}(\omega').$$

The information drift of $(\mathcal{G}_t)$ relative to $(\mathcal{F}_t)$ is already determined by the trace of $(\beta_t)$. For the corresponding analogue of Theorem 2.6, we shall give a more explicit statement.

THEOREM 2.12. *The process*

$$\alpha_t(\omega) = \beta_t(\omega, \omega)$$

*is the information drift of $(\mathcal{G}_t)$ relative to $(\mathcal{F}_t)$.*

PROOF. Suppose $T$ to be a stopping time such that $M^T$ is a martingale. For $0 \leq s < t$, $A \in \mathcal{H}_s^0$ and $B \in \mathcal{F}_s^0$, we have to show

$$E[\mathbb{1}_A \mathbb{1}_B (M_t^T - M_s^T)] = E\left[\mathbb{1}_A \mathbb{1}_B \int_s^t \beta_u(\omega, \omega) \, d\langle M, M \rangle_u^T \right].$$

For simplicity, assume $M^T = M$, and observe, like in the proof of Theorem 2.6,

$$E[\mathbb{1}_A \mathbb{1}_B (M_t - M_s)] = E[\mathbb{1}_B P_t(\cdot, A)(M_t - M_s)]$$
$$= E\left[\mathbb{1}_A(\omega) \mathbb{1}_B(\omega) \int_s^t \beta_u(\omega, \omega) \, d\langle M, M \rangle_u \right]. \qquad \square$$

EXAMPLE 2.13. Let $(W_t)$ be the standard Wiener process and $(\mathcal{F}_t^0)$ the filtration generated by $(W_t)$. Moreover, let $(Y_t)$ be a Gaussian process independent of $\mathcal{F}_1$ such that, for each pair $s,t$ with $0 \leq s < t$, the difference $Y_t - Y_s$ is independent of $Y_t$. We denote by $w_t$ the variance of $Y_t$.

We enlarge our filtration by
$$\mathcal{H}_t^0 = \sigma(W_1 + Y_s : 0 \leq s \leq t) = \sigma(W_1 + Y_t) \vee \sigma(Y_t - Y_s : 0 \leq s \leq t),$$
and put $\mathcal{G}_t^0 = \mathcal{F}_t^0 \vee \mathcal{H}_t^0$, $0 \leq t \leq 1$. Now observe that, for all $C \in \sigma(Y_t - Y_s : 0 \leq s \leq t)$ and Borel sets $B \in \mathcal{B}(\mathbb{R})$, we have

$$P_t(\cdot, \{W_1 + Y_t \in B\} \cap C) = P(C) \int \mathbb{1}_B(x + W_1 - W_t + Y_t) \, dP \Big|_{x=W_t}$$
$$= P(C) \int \mathbb{1}_B(y + x) \phi_{1-t+w_t}(y) \, dy \Big|_{x=W_t}$$
$$= P(C) \int_B \phi_{1-t+w_t}(y - W_t) \, dy, \qquad 0 \leq t < 1,$$

where
$$\phi_v(y) = \frac{1}{(2\pi v)^{1/2}} e^{-y^2/(2v)}.$$

Now observe that $f(x,t) = P(C) \int_B \phi_{1-t+w_t}(y-x)\,dy$ is differentiable in $x$ and satisfies
$$\frac{\partial}{\partial x} f(x,t) = P(C) \int_B \frac{y-x}{1-t+w_t} \phi_{1-t+w_t}(y-x)\,dy$$
for all $0 \leq t < 1$ and $x \in \mathbb{R}$. By Itô's formula,
$$P_t(\cdot, \{W_1 + Y_t \in B\} \cap C) = f(0,0) + \int_0^t \frac{\partial}{\partial x} f(W_s, s)\,dW_s + A_t, \qquad 0 \leq t < 1,$$
where $A$ is a process of bounded variation. Note that $A$ is also a martingale and, thus, $A = 0$. Hence,
$$k_t(\cdot, \{W_1 + Y_t \in B\} \cap C)$$
$$= P(C) \int_B \frac{y - W_t}{1-t+w_t} \phi_{1-t+w_t}(y - W_t)\,dy$$
$$= P(C) \int \mathbb{1}_B(y+x) \frac{y+x-x}{1-t+w_t} \phi_{1-t+w_t}(y)\,dy \Big|_{x=W_t(\omega)}$$
$$= \int_{\{W_1+Y_t\in B\}\cap C} \frac{W_1(\omega') + Y_t(\omega') - W_t(\omega)}{1-t+w_t}\,dP_t(\omega, d\omega').$$

As a consequence,
$$\beta_t(\omega, \omega') = \frac{k_t(\omega, d\omega')}{P_t(\omega, d\omega')}\Big|_{\mathcal{H}_t^0} = \frac{W_1(\omega') + Y_t(\omega') - W_t(\omega)}{1-t+w_t},$$

and by Theorem 2.12,

$$W_t - \int_0^t \frac{W_1 + Y_s - W_s}{1 - s + w_s} \, ds, \qquad 0 \leq t < 1,$$

is a martingale relative to $(\mathcal{G}_t)$.

Similar examples can be found in [8] where the information drifts are derived in a completely different way though.

EXAMPLE 2.14. Let $(W_t)$ be the standard Wiener process and $(\mathcal{F}_t)$ the Wiener filtration. We use the abbreviation $W_t^* = \sup_{0 \leq s \leq t} W_s$ and consider the filtration enlarged by the random variable $G = \mathbb{1}_{[0,c]}(W_1^*)$, $c > 0$. Again, we want to apply Theorem 2.12 in order to obtain the information drift of $\mathcal{G}_t = \mathcal{F}_t \vee \sigma(G)$. To this end, let $Z_t = \sup_{t \leq r \leq 1}(W_r - W_t)$ and denote by $p_t$ the density of $Z_t$, $0 \leq t < 1$. Now,

$$P_t(\cdot, G = 1) = P(W_t^* \vee W_t + Z_t \leq c | \mathcal{F}_t)$$
$$= \int \mathbb{1}_{[0,c]}(y \vee x + Z_t) \, dP \bigg|_{x = W_t, y = W_t^*}$$
$$= \mathbb{1}_{[0,c]}(y) \int_0^{c-x} p_t(z) \, dz \bigg|_{x = W_t, y = W_t^*},$$

for all $0 \leq t < 1$. Note that $F(x, y, t) = \mathbb{1}_{[0,c]}(y) \int_0^{c-x} p_t(z) \, dz$ is differentiable in $x$ for all $0 \leq t < 1$ and $x \in \mathbb{R}$, and by Itô's formula,

$$P_t(\cdot, G = 1) = F(0, 0, 0) + \int_0^t \frac{\partial}{\partial x} F(W_s, W_s^*, s) \, dW_s + A_t, \qquad 0 \leq t < 1,$$

where $A$ is a process of bounded variation. Hence,

$$k_t(\cdot, G = 1) = \frac{\partial}{\partial x} F(W_t, W_t^*, t), \qquad 0 \leq t < 1.$$

Similarly, we have

$$P_t(\cdot, G = 0) = H(W_t, W_t^*, t), \qquad 0 \leq t < 1,$$

and

$$k_t(\cdot, G = 0) = \frac{\partial}{\partial x} H(W_t, W_t^*, t), \qquad 0 \leq t < 1,$$

where

$$H(x, y, t) = \mathbb{1}_{(c, \infty)}(y) + \mathbb{1}_{[0,c]}(y) \int_{c-x}^\infty p_t(z) \, dz.$$

As a consequence,

$$\begin{aligned}\beta_t(\omega,\omega') &= \frac{k_t(\omega,d\omega')}{P_t(\omega,d\omega')}\bigg|_{\sigma(G)} \\ &= \mathbb{1}_{\{1\}}(G(\omega'))\frac{\partial}{\partial x}\log F(W_t(\omega),W_t^*(\omega'),t) \\ &\quad + \mathbb{1}_{\{0\}}(G(\omega'))\frac{\partial}{\partial x}\log H(W_t(\omega),W_t^*(\omega'),t), \qquad 0\le t<1.\end{aligned}$$

**3. Monotone convergence of information drifts.** In the preceding section we established a general relationship between the information drift and the regular conditional probabilities of filtrations. In this framework the knowledge of the better-informed agent is described by a general enlarged filtration $(\mathcal{G}_t)$ of $(\mathcal{F}_t)$. We shall now consider the question whether this situation may be well approximated by "stepwise initial" enlargements, for which we take $\mathcal{F}_t \vee \mathcal{G}_{t_i-}$ for $t \in [t_i, t_{i+1})$, if the family $(t_i)_{0 \le i \le n}$ is a partition of $\mathbb{R}_+$. One particularly important question in this context concerns the behavior of the information drifts along such a sequence of discretized enlargements. Of course, we expect some convergence of the drifts. We shall establish this fact rigorously in the following section. In the present section, we shall prepare the treatment of this problem by solving a somewhat more general problem. Let $(\mathcal{G}_t^n)_{n \in \mathbb{N}}$ be an increasing sequence of finite utility filtrations and $\sup_n u_{\mathcal{G}^n}(x)$ be finite. We will show that the smallest filtration containing every $(\mathcal{G}_t^n)$ is then also a finite utility filtration.

Since we will not deal with regular conditional probabilities in this section, it is not necessary to require our probability space $(\Omega, \mathcal{F})$ to be standard.

We use the terminology of Revuz and Yor [24]: $H^2(\mathcal{F}_t)$ denotes the set of $L^2$-bounded continuous $(\mathcal{F}_t)$-martingales, that is, the space of continuous $(\mathcal{F}_t, P)$-martingales $M$ such that

$$\sup_{t \ge 0} E(M_t^2) < \infty.$$

We need the following characterization of $H^2(\mathcal{F}_t)$.

LEMMA 3.1 (Proposition 1.23 in [24]). *A continuous $(\mathcal{F}_t)$-local martingale belongs to $H^2(\mathcal{F}_t)$ if and only if the following two conditions hold:*

(i) $E(M_0^2) < \infty$,
(ii) $E(\langle M, M \rangle_\infty) < \infty$.

The properties (i) and (ii) are independent of the filtration considered. This is due to the fact that the quadratic variation of $M$ does not change under a new filtration $(\mathcal{G}_t)$ for which $M$ is still a semimartingale. We therefore have the following:

LEMMA 3.2. *Suppose $M \in H^2(\mathcal{F}_t)$. Let $(\mathcal{G}_t)$ be a filtration such that $M$ is still a $(\mathcal{G}_t)$-semimartingale. If*

$$M = \widetilde{M} + A$$

*is a Doob–Meyer decomposition with respect to $(\mathcal{G}_t)$ with $A_0 = 0$, then $\widetilde{M}$ belongs to $H^2(\mathcal{G}_t)$.*

PROOF. Notice that $\widetilde{M}_0 = M_0$ and $\langle M, M \rangle = \langle \widetilde{M}, \widetilde{M} \rangle$. The claim follows now by applying Lemma 3.1 twice. □

Now let $M$ be a continuous $(\mathcal{F}_t)$-local martingale and $(\mathcal{G}_t^n)_{n \geq 1}$ an increasing sequence of filtrations, that is, for all $t \geq 0$, we have

$$\mathcal{F}_t \subset \mathcal{G}_t^1 \subset \cdots \subset \mathcal{G}_t^n \subset \mathcal{G}_t^{n+1} \subset \cdots.$$

We assume that, for all $n \geq 1$, the process $M$ is a $(\mathcal{G}_t^n)$-semimartingale with Doob–Meyer decomposition of the form

$$M = M^n + \int_0^\cdot \mu_s^n \, d\langle M, M \rangle_s,$$

where $\mu^n$ is $(\mathcal{G}_t^n)$-predictable. We then have the following asymptotic property.

LEMMA 3.3. *If the processes $(\mu^n)_{n \in \mathbb{N}}$ converge to some $\mu$ in $L^2(M)$, then*

$$M - \int_0^\cdot \mu_s \, d\langle M, M \rangle_s$$

*is a local martingale with respect to $\mathcal{G}_t = \bigvee_{n \geq 1} \mathcal{G}_t^n, t \geq 0$.*

PROOF. Suppose the stopping time $\tau$ reduces $M$ such that $M^\tau$ is a bounded martingale. Note that Lemma 3.2 implies that the stopped processes $(M^n)^\tau$ are $(\mathcal{G}_t)$-martingales.

For simplicity, we assume $M^\tau = M$. Choose a constant $C > 0$ such that

$$|M| \leq C \quad \text{and} \quad E \int_0^\infty (\mu_s^n)^2 \, d\langle M, M \rangle_s \leq C^2 \qquad \text{for all } n \geq 1.$$

Now let $\varepsilon > 0$, $0 \leq s < t$ and $A \in \mathcal{G}_s$. It suffices to show

$$\left| E[\mathbb{1}_A(M_t - M_s)] - E\left[\mathbb{1}_A \int_s^t \mu_s \, d\langle M, M \rangle_s\right] \right| \leq \varepsilon.$$

We start by choosing $n_0$ such that

$$\|\mu^n - \mu\|_{L^2(M)} \leq \frac{\varepsilon}{4\sqrt{E(\langle M, M \rangle)_\infty}}$$

for all $n \geq n_0$.

Note that $\bigcup_{n \geq n_0} \mathcal{G}_s^n$ is an algebra generating the $\sigma$-algebra $\bigvee_{n \geq 1} \mathcal{G}_s^n = \mathcal{G}_s = \bigvee_{n \geq n_0} \mathcal{G}_s^n$. Hence, we can find a sequence $(A_i)_{i \in \mathbb{N}}$ of sets in $\bigcup_{n \geq n_0} \mathcal{G}_s^n$ such that $P(A \triangle A_i) \to 0$. A subsequence of $(\mathbb{1}_{A_i})_{i \in \mathbb{N}}$ converges to $\mathbb{1}_A$ almost surely and, therefore, we may choose $n \geq n_0$ and $\tilde{A} \in \mathcal{G}_s^n$ satisfying $P(\tilde{A} \triangle A) \leq (\frac{\varepsilon}{4C})^2$ and

$$\left( E \int_s^t (\mathbb{1}_A - \mathbb{1}_{\tilde{A}})^2 \, d\langle M, M \rangle \right)^{1/2} \leq \frac{\varepsilon}{4C}.$$

Hence, we have

$$|E[\mathbb{1}_A(M_t - M_s)] - E[\mathbb{1}_{\tilde{A}}(M_t - M_s)]| \leq |E[(\mathbb{1}_{\tilde{A}} - \mathbb{1}_A)(M_t - M_s)]|$$
$$\leq P(\tilde{A} \triangle A)^{1/2} (E(M_t - M_s)^2)^{1/2}$$
$$\leq \frac{\varepsilon}{2}.$$

By applying the Kunita–Watanabe inequality, we get, for $n \geq n_0$,

$$\left| E\left[ \mathbb{1}_A \int_s^t \mu_u \, d\langle M, M \rangle_u \right] - E\left[ \mathbb{1}_{\tilde{A}} \int_s^t \mu_u^n \, d\langle M, M \rangle_u \right] \right|$$
$$\leq \left| E\left[ \mathbb{1}_A \int_s^t (\mu_u - \mu_u^n) \, d\langle M, M \rangle_u \right] \right| + \left| E\left[ (\mathbb{1}_A - \mathbb{1}_{\tilde{A}}) \int_s^t \mu_u^n \, d\langle M, M \rangle_u \right] \right|$$
$$\leq \left( E \int_s^t \mathbb{1}_A \, d\langle M, M \rangle \right)^{1/2} \left( E \int_0^\infty (\mu_s - \mu_s^n)^2 \, d\langle M, M \rangle_s \right)^{1/2}$$
$$+ \left( E \int_s^t (\mathbb{1}_A - \mathbb{1}_{\tilde{A}})^2 \, d\langle M, M \rangle \right)^{1/2} \left( E \int_0^\infty (\mu_s^n)^2 \, d\langle M, M \rangle_s \right)^{1/2}$$
$$\leq (E\langle M, M \rangle)^{1/2} \|\mu - \mu^n\|_{L^2(M)} + \frac{\varepsilon}{4}$$
$$\leq \frac{\varepsilon}{2},$$

and thus,

$$\left| E[\mathbb{1}_A(M_t - M_s)] - E\left[ \mathbb{1}_A \int_s^t \mu_u \, d\langle M, M \rangle_u \right] \right|$$
$$\leq |E[\mathbb{1}_A(M_t - M_s)] - E[\mathbb{1}_{\tilde{A}}(M_t - M_s)]|$$
$$+ \left| E[\mathbb{1}_{\tilde{A}}(M_t - M_s)] - E\left[ \mathbb{1}_{\tilde{A}} \int_s^t \mu_u^n \, d\langle M, M \rangle_u \right] \right|$$
$$+ \left| E\left[ \mathbb{1}_{\tilde{A}} \int_s^t \mu_u^n \, d\langle M, M \rangle_u \right] - E\left[ \mathbb{1}_A \int_s^t \mu_u \, d\langle M, M \rangle_u \right] \right|$$
$$\leq \frac{\varepsilon}{2} + 0 + \frac{\varepsilon}{2} = \varepsilon.$$

□

We are now in a position to prove the main result of the section.

THEOREM 3.4. *If* $\sup_{n\geq 1} \|\mu^n\|^2_{L^2(M)} < \infty$, *then* $(\mu^n)$ *converges in* $L^2(M)$ *to a process* $\mu$. *Moreover*,

$$M - \int_0^\cdot \mu \, d\langle M, M \rangle$$

*is a local martingale with respect to* $\mathcal{G}_t = \bigvee_{n\geq 1} \mathcal{G}_t^n, t \geq 0$.

PROOF. Set $c = \sup_{n\geq 1} \|\mu^n\|^2_{L^2(M)}$. Let $m \geq n \geq 1$, and note that $\mu^m - \mu^n$ is the information drift of $(\mathcal{G}_t^m)$ relative to $(\mathcal{G}_t^n)$. Therefore, property (iii) in Lemma 1.6 implies

$$\|\mu^m\|^2_{L^2(M)} = \|\mu^n\|^2_{L^2(M)} + \|\mu^m - \mu^n\|^2_{L^2(M)}.$$

Thus, $c = \lim_{n\to\infty} \|\mu^n\|^2_{L^2(M)}$ and

$$\|\mu^m - \mu^n\|^2_{L^2(M)} = \|\mu^m\|^2_{L^2(M)} - \|\mu^n\|^2_{L^2(M)} \leq c - \|\mu^n\|^2_{L^2(M)} \to 0$$

as $n \to \infty$. Therefore, $\{\mu^n\}_{n\geq 1}$ is a Cauchy sequence in $L^2(M)$. By completeness of $L^2(M)$, there exists a unique $(\mathcal{G}_t)$-predictable process $\mu^0 \in L^2(M)$ such that $\lim_{n\to\infty} \mu^n = \mu^0$ in $L^2(M)$. By Lemma 3.3, the process $M - \int \mu^0 \, d\langle M, M \rangle$ is a $(\mathcal{G}_t)$-local martingale. □

**4. Continuous and initial enlargements.** In this section we relate general enlargements of filtrations to "initial enlargements" along discrete partitions of $[0, T]$, for finite horizon $T$. The knowledge of the insider is modeled by an arbitrary filtration $(\mathcal{G}_t)_{t\in[0,T]}$, satisfying $\mathcal{G}_t \supset \mathcal{F}_t, 0 \leq t \leq T$. For $s \in [0, T]$, we set

$$\mathcal{G}_t^s = \begin{cases} \mathcal{F}_t, & t < s, \\ \mathcal{F}_t \vee \mathcal{G}_{s-}, & t \geq s. \end{cases}$$

Again, the analysis of this section does not require our probability space $(\Omega, \mathcal{F})$ to be standard.

REMARK 4.1. In the case where the $\sigma$-field $\mathcal{G}_{s-}, s \in [0, T]$, is generated by a countable number of events, say, $(A_n)_{n\in\mathbb{N}}$, the enlarged filtration $\mathcal{G}_t^s$ can be viewed as initial enlargement at time $s$ in the classical sense. In that case $\mathcal{G}_{s-} = \sigma(\mathbb{1}_{A_n} : n \in \mathbb{N})$ and one has, for $t \in [0, T]$,

$$\mathcal{F}_t \vee \mathcal{G}_{s-} = \mathcal{F}_t \vee \sigma(\mathbb{1}_{A_n} : n \in \mathbb{N}).$$

The set $\{0, 1\}^\mathbb{N}$ can be endowed with a metric so that it becomes a Polish space with corresponding Borel$-\sigma$-field $\mathcal{B}(\{0, 1\})^{\otimes \mathbb{N}}$. Hence, the filtration $(\mathcal{G}_t^s)$ can be seen as initial enlargement at time $s$ induced by the random variable $G:\Omega \to \{0,1\}^{\mathbb{N}}, \omega \mapsto (\mathbb{1}_{A_n}(\omega))_{n\in\mathbb{N}}$. In particular, the standard theory of initial enlargements is applicable. See [21].

In the following, we assume that $(\mathcal{G}_t^s)$ is for arbitrary $s \in [0,T]$ a finite utility filtration. For $0 \le s \le t \le T$, we denote

$$\pi_0([0,s) \times (t,T]) = F(s,t) = \tfrac{1}{2} E \int_t^T (\mu_r^s)^2 \, d\langle M, M \rangle_r,$$

where $\mu^s$ is a $(\mathcal{G}_t^s)$-information drift. So far $\pi_0$ is defined only on the set $J = \{[0,s) \times (t,T] : s \le t\}$. As the next lemma shows, $\pi_0$ can be extended to a measure on the Borel sets of $D = \{(s,t) \in \mathbb{R}^2 : 0 \le s < t \le T\}$.

LEMMA 4.2. *There exists a unique measure $\pi$ on the Borel sets $\mathcal{B}(D)$ of $D$ satisfying $\pi|_J = \pi_0$.*

PROOF. Uniqueness is an immediate consequence of the measure extension theorem. In order to show the existence of an extension, it satisfies to verify the following property which essentially amounts to countable additivity on a generating semiring (see [13], Chapter II, Satz 3.8): For any $(s,t) \in \overline{D}$ and any sequence $(s_n,t_n)_{n\in\mathbb{N}}$ in $\overline{D}$ with $s_n \le s$, $t_n \ge t$ and $\lim_{n\to\infty}(s_n,t_n) = (s,t)$, we have $\lim_{n\to\infty} F(s_n,t_n) = F(s,t)$. Moreover, $F(s_n,t_n) \le F(s,t) < \infty$.

Let $s_n$, $t_n$, $s$ and $t$ as above. Without loss of generality, we assume that $(s_n)$ is monotonically increasing. For $u \in [t,T]$, we consider the filtrations $(\mathcal{G}_r^{s_n})_{r\in[u,T]}$, $n \in \mathbb{N}$, over the time interval $[u,T]$. The filtrations are monotonically increasing with $\bigvee_{n\in\mathbb{N}} \mathcal{G}_r^{s_n} = \mathcal{G}_r^s$, $r \in [u,T]$. Since $(\mu_r^{s_n})_{r\in[u,T]}$ are $(\mathcal{G}_r^{s_n})$-information drifts, it follows (by Lemma 1.6) that

$$E \int_u^T (\mu^s - \mu^{s_n}) \mu^{s_n} \, d\langle M, M \rangle = 0.$$

In particular,

$$E \int_u^T (\mu^{s_n})^2 \, d\langle M, M \rangle \le E \int_t^T (\mu^s)^2 \, d\langle M, M \rangle < \infty.$$

By Theorem 3.4, the processes $(\mu_r^{s_n})_{r\in[u,T]}$ converge to the information drift $(\mu_r^s)_{r\in[u,T]}$ in $L^2(M;[u,T])$. Therefore, for any $u \in (t,T]$,

$$\liminf_{n\to\infty} E \int_{t_n}^T (\mu^{s_n})^2 \, d\langle M, M \rangle \ge E \int_u^T (\mu^s)^2 \, d\langle M, M \rangle.$$

Due to the continuity of $M$, the right-hand side of the previous equation tends to $E \int_t^T (\mu^s)^2 \, d\langle M, M \rangle$ as $u \downarrow t$. Consequently, we obtain $\lim_{n\to\infty} F(s_n,t_n) = F(s,t)$. □

The measure $\pi$ describes the utility increase by additional information. As will be shown below, $\pi(D)$ is finite if and only if $(\mathcal{G}_t)$ is a finite utility filtration.

We now approximate the general filtration $(\mathcal{G}_t)$ by filtrations that can be seen as successive initial enlargements. Let $\Delta : 0 = s_0 \leq \cdots \leq s_n = T$, $n \in \mathbb{N}$, be a partition of the interval $[0,T]$. We let, for $r \in [s_i, s_{i+1})$, $i = 0, \ldots, n-1$,
$$\mathcal{G}_r^\Delta = \mathcal{G}_{s_i-} \vee \mathcal{F}_r.$$

PROPOSITION 4.3. *For $i = 0, \ldots, n-1$, let $\mu^{s_i}$ be a $(\mathcal{G}_r^{s_i})$-information drift and set $\mu_r^\Delta := \mu_r^{s_i}$ for $r \in [s_i, s_{i+1})$. Then $\mu^\Delta$ is a $\mathcal{G}_t^\Delta$-information drift. Moreover,*
$$\tfrac{1}{2} \int_0^T (\mu_r^\Delta)^2 \, d\langle M, M \rangle_r = \pi(D^\Delta),$$
*where $D^\Delta := \{(s,t) \in D : \exists i \in \{0, \ldots, n-1\} \text{ with } s < s_i \text{ and } t > s_i\}$.*

PROOF. It is straightforward to verify that $\mu^\Delta$ is an information drift for $(\mathcal{G}_t^\Delta)$. Moreover,
$$\tfrac{1}{2} E \int_0^T (\mu_r^\Delta)^2 \, d\langle M, M \rangle_r$$
$$= \tfrac{1}{2} \sum_{i=0}^{n-1} E \int_{s_i}^{s_{i+1}} (\mu_r^{s_i})^2 \, d\langle M, M \rangle_r$$
$$= \tfrac{1}{2} \sum_{i=0}^{n-1} \left( E \int_{s_i}^T (\mu_r^{s_i})^2 \, d\langle M, M \rangle_r - E \int_{s_{i+1}}^T (\mu_r^{s_i})^2 \, d\langle M, M \rangle_r \right)$$
$$= \sum_{i=0}^{n-1} \pi([0, s_i) \times (s_i, s_{i+1}]) = \pi(D^\Delta). \qquad \square$$

We can now state the main theorem of this section.

THEOREM 4.4. *Let $\Delta_n$, $n \in \mathbb{N}$, be a sequence of partitions of the interval $[0,T]$, the mesh of which tends to $0$. If $\pi(D)$ is finite, then the information drifts $\mu^{\Delta_n}$ converge in $L^2(M)$ to a $(\mathcal{G}_t)$-information drift $\mu$. Moreover, the utility gain of the insider satisfies*
$$u(\mathcal{G}_t, x) - u(\mathcal{F}_t, x) = \tfrac{1}{2} E \int_0^T \mu^2 \, d\langle M, M \rangle = \pi(D).$$
*If $\pi(D)$ is infinite, then so is the utility gain of the insider.*

The proof of the theorem is based on the following proposition.

PROPOSITION 4.5. *If $\pi(D) < \infty$, then there exists a $(\mathcal{G}_t)$-information drift $\mu$. Moreover,*

$$\tfrac{1}{2}\|\mu\|_{L^2(M)}^2 = \pi(D).$$

PROOF. Let $\Delta_n$, $n \in \mathbb{N}$, be as in the above theorem with the additional assumption that $\Delta_{n+1}$ is a refinement of $\Delta_n$ for all $n \in \mathbb{N}$. Then one has $\mathcal{G}_t^{\Delta_n} \subset \mathcal{G}_t^{\Delta_{n+1}}$ for any $t \in [0, T]$. By Proposition 4.3, $\tfrac{1}{2}\|\mu^{\Delta_n}\|_{L^2(M)}^2 = \pi(D^{\Delta_n}) \leq \pi(D)$. Due to Theorem 3.4, the information drifts $\mu^{\Delta_n}$ converge to a $(\bigvee_{n \in \mathbb{N}} \mathcal{G}_t^{\Delta_n}) = (\mathcal{G}_{t-})$-information drift $\mu$ in $L^2(M)$. Using monotone convergence, we obtain that

$$\pi(D) = \lim_{n \to \infty} \pi(D^{\Delta_n}) = \lim_{n \to \infty} \tfrac{1}{2}\|\mu^{\Delta_n}\|_{L^2(M)}^2 = \tfrac{1}{2}\|\mu\|_{L^2(M)}^2.$$

Since every cadlag $(\mathcal{G}_{t-})$-martingale is as well a $(\mathcal{G}_t)$-martingale, $\mu$ is a $(\mathcal{G}_t)$-information drift. □

PROOF OF THEOREM 4.4. Assume that $\pi(D)$ is finite. Since the mesh of the partitions $\Delta_n$ tends to zero, one has $\lim_{n \to \infty} \mathbb{1}_{D^{\Delta_n}}(x) = 1$ for all $x \in D$. Consequently, the dominated convergence theorem yields

$$\lim_{n \to \infty} \pi(D^{\Delta_n}) = \pi(D). \tag{11}$$

We established the existence of a $(\mathcal{G}_t)$-information drift $\mu$ in Proposition 4.5. Recall that, by Lemma 1.6, the processes $\mu^{\Delta_n}$ and $\mu - \mu^{\Delta_n}$ are orthogonal in $L^2(M)$. Consequently,

$$\|\mu - \mu^{\Delta_n}\|_{L^2(M)}^2 = \|\mu\|_{L^2(M)}^2 - \|\mu^{\Delta_n}\|_{L^2(M)}^2.$$

Due to (11), the right-hand side of the previous equation converges to 0. Hence, $\mu^{\Delta_n}$ converges to $\mu$ in $L^2(M)$. The remaining statements are consequences of Propositions 4.5 and 1.4. □

**5. Additional utility and entropy of filtrations.** In this section we consider the link between the additional expected logarithmic utility of a better informed agent and the entropy of the additional information he possesses. The additional utility was first expressed in terms of a relative entropy in [23], page 1103, for a particular example. More generally, [3] discussed the link between the absolute entropy of a random variable describing initially available additional information and the utility increment of better informed agents. Here we shall see that the expected logarithmic utility increment is given by an integral version of relative entropies of the $\sigma$-algebras of the filtration. This notion can best be understood as the limit of discrete entropy sums along a sequence of partitions of the trading interval as the mesh goes to 0. Alternatively, we are able to give an interpretation of the utility increment by *Shannon information* differences between the filtrations of the agents. In particular, we shall see that these differences are independent of any local martingales the filtrations may carry.

Suppose that the assumptions of Section 2 are satisfied. Moreover, we assume that $M$ is a continuous local martingale satisfying the (PRP) relative to $(\mathcal{F}_t)$, which simply means that $L^A = 0$. Equation (7) simplifies to

$$P_t(\cdot, A) = P(A) + \int_0^t k_s(\cdot, A) \, dM_s,$$

where $k(\cdot, A)$ is as in Section 2. Let again $(\mathcal{G}_t^0)$ be a filtration satisfying $\mathcal{F}_t^0 \subset \mathcal{G}_t^0$ and being generated by countably many sets. To simplify notation, we assume the filtration $(\mathcal{G}_t^0)$ to be left-continuous. Let $(\mathcal{G}_t)$ be the smallest completed and right-continuous filtration containing $(\mathcal{G}_t^0)$. In the following, we assume that $(\mathcal{G}_t)$ is a finite utility filtration and denote by $\mu$ its predictable information drift, that is,

$$\widetilde{M} = M - \int_0^\cdot \mu_t \, d\langle M, M \rangle_t$$

is a $(\mathcal{G}_t)$-local martingale. Recall that, by Theorem 2.10, we may assume that $k_t(\omega, \cdot)$ is a signed measure. For a fixed $r > 0$, we define $\mu^r$ as the information drift of the initially enlarged filtration $(\mathcal{G}_t^r)$, defined as in the beginning of the preceding chapter. For stating the main result we need the following lemma.

LEMMA 5.1. *Let $0 \leq s < t$ and $(\mathcal{P}^m)_{m \geq 0}$ be an increasing sequence of finite partitions such that $\sigma(\mathcal{P}^m : m \geq 0) = \bar{\mathcal{G}}_s^0$. Then*

$$\lim_m E \int_s^t \sum_{A \in \mathcal{P}^m} \left(\frac{k_u}{P_u}\right)^2 (\cdot, A) \mathbb{1}_A \, d\langle M, M \rangle_u = E \int_s^t (\mu_u^s)^2 \, d\langle M, M \rangle_u$$

*and*

$$\lim_m E \int_s^t \sum_{A \in \mathcal{P}^m} \frac{k_u}{P_u}(\cdot, A) \mathbb{1}_A \mu_u^s \, d\langle M, M \rangle_u = E \int_s^t (\mu_u^s)^2 \, d\langle M, M \rangle_u.$$

PROOF. By Lemma 2.10, the process

$$Y_u^m(\omega, \omega') = \sum_{A \in \mathcal{P}^m} \frac{k_u}{P_u}(\omega, A) \mathbb{1}_A(\omega'), \qquad m \geq 1,$$

is a $L^2$-bounded martingale for $P_M$-a.a. $(\omega, u) \in \Omega \times [s, t]$. Hence, $(Y^m)$ converges $P_M$-a.s. to the density

$$\gamma_u = \left.\frac{k_u(\cdot, d\omega')}{P_u(\cdot, d\omega')}\right|_{\mathcal{G}_s^0}.$$

By Theorem 2.6, we have

$$\gamma_u(\omega, \omega') = \mu_u^s(\omega)$$

$P_M$-a.s. on $\Omega \times [s,t]$ and, hence, the first result. In a similar way, one can prove the second statement. □

We next discuss the important concept of the additional information of a $\sigma$-field relative to a filtration.

DEFINITION 5.2. Let $\mathcal{A}$ be a sub-$\sigma$-algebra of $\mathcal{F}$ and $R, Q$ two probability measures on $\mathcal{F}$. Then we define the relative entropy of $R$ with respect to $Q$ on the $\sigma$-field $\mathcal{A}$ by

$$\mathcal{H}_\mathcal{A}(R\|Q) = \begin{cases} \int \log \left.\frac{dR}{dQ}\right|_\mathcal{A} dR, & \text{if } R \ll Q, \\ \infty, & \text{else.} \end{cases}$$

Moreover, the *additional information* of $\mathcal{A}$ relative to the filtration $(\mathcal{F}_r)$ on $[s,t]$ $(0 \leq s < t \leq T)$ is defined by

$$H_\mathcal{A}(s,t) = \int \mathcal{H}_\mathcal{A}(P_t(\omega, \cdot) \| P_s(\omega, \cdot)) \, dP(\omega).$$

The following lemma establishes the basic link between the entropy of a filtration enlargement and additional logarithmic utility of a trader possessing this information advantage.

LEMMA 5.3. *For $0 \leq s < t$, we have*

$$H_{\mathcal{G}_s^0}(s,t) = \tfrac{1}{2} E \int_s^t (\mu_u^s)^2 \, d\langle M, M \rangle_u.$$

PROOF. Let $(\mathcal{P}^m)_{m \geq 0}$ be an increasing sequence of finite partitions such that $\sigma(\mathcal{P}^m : m \geq 0) = \mathcal{G}_s^0$. Recall that, by (9)

$$\sum_{A \in \mathcal{P}^m} \left[ \mathbb{1}_A \log P_s(\cdot, A) - \mathbb{1}_A \log P_t(\cdot, A) \right]$$

$$= \sum_{A \in \mathcal{P}^m} \left[ -\int_s^t \frac{k_u}{P_u}(\cdot, A) \mathbb{1}_A \, d\widetilde{M}_u - \int_s^t \frac{k_u}{P_u}(\cdot, A) \mathbb{1}_A \mu_u \, d\langle M, M \rangle_u \right.$$

$$\left. + \frac{1}{2} \int_s^t \left(\frac{k_u}{P_u}\right)^2 (\cdot, A) \mathbb{1}_A \, d\langle M, M \rangle_u \right].$$

Since $\widetilde{M}$ is a local martingale, we obtain by stopping and taking limits if necessary

$$E \sum_{A \in \mathcal{P}^m} P_s(\cdot, A) \log \frac{P_t(\cdot, A)}{P_s(\cdot, A)}$$

$$= E \sum_{A \in \mathcal{P}^m} \int_s^t \frac{k_u}{P_u}(\cdot, A) \mathbb{1}_A \mu_u \, d\langle M, M \rangle_u - \frac{1}{2} \int_s^t \left(\frac{k_u}{P_u}\right)^2 (\cdot, A) \mathbb{1}_A \, d\langle M, M \rangle_u.$$

Note that in the previous line $\mu$ may be replaced by $\mu^s$, because $(\mu - \mu^s)$ is orthogonal to $L^2(M)(\mathcal{G}^s)$ [see property (iv) in Lemma 1.6]. Applying Lemma 5.1 yields

$$\lim_m H_{\mathcal{P}^m}(s,t) = \tfrac{1}{2} E \int_s^t (\mu_u^s)^2 \, d\langle M, M \rangle_u.$$

Fatou's lemma implies

$$\liminf_m H_{\mathcal{P}^m}(s,t) \geq H_{\mathcal{G}_s^0}(s,t).$$

On the other hand, we have $H_{\mathcal{P}^m}(s,t) \leq H_{\mathcal{G}_s^0}(s,t)$, since $\mathcal{P}^m \subset \mathcal{G}_s^0$ and, thus,

$$\lim_m H_{\mathcal{P}^m}(s,t) = H_{\mathcal{G}_s^0}(s,t),$$

which completes the proof. $\square$

Let us now return to the stepwise approximation of a filtration enlargement along a sequence of partitions of the trading interval by "initial enlargements," and define their respective information increment.

DEFINITION 5.4. Let $\Delta : 0 = s_0 \leq \cdots \leq s_n = T$, $n \in \mathbb{N}$, be a partition of the interval $[0, T]$ and let $\mu^\Delta$ be the information drift of $(\mathcal{G}_r^\Delta)$. The additional information of $(\mathcal{G}_r^\Delta)$ relative to $(\mathcal{F}_r)$ is defined as

$$H_\Delta = \sum_{i=0}^{n-1} H_{\mathcal{G}_{s_i}^0}(s_i, s_{i+1}).$$

THEOREM 5.5. *We have*

$$\lim_{|\Delta| \to 0} H_\Delta = \tfrac{1}{2} E \int_0^T \mu_u^2 \, d\langle M, M \rangle_u.$$

PROOF. This follows from Theorem 4.4 and Lemma 5.3. $\square$

EXAMPLE 5.6. Let $\mathcal{G}_t^0 = \mathcal{F}_t^0 \vee \sigma(\mathcal{P})$, where $\mathcal{P}$ is a finite partition in $\mathcal{F}_T$. Then $\mu^0 = \mu$ and by Lemma 5.3,

$$H_{\mathcal{G}_0^0}(0,T) = \tfrac{1}{2} E \int_0^T \mu_u^2 \, d\langle M, M \rangle_u.$$

If $\mathcal{F}_0$ is trivial, then

$$H_{\mathcal{G}_0^0}(0,T) = -\sum_{A \in \mathcal{P}} P(A) \log P(A),$$

which is the absolute entropy of the partition $\mathcal{P}$. Thus, the additional logarithmic utility of an agent with information $(\mathcal{G}_t)$ is equal to the entropy of $\mathcal{P}$. This example shows that there is a link between logarithmic utility and the so-called Shannon information.

DEFINITION 5.7. Let $X$ and $Y$ be two random variables in some measurable spaces. The *mutual information* (or Shannon information) between $X$ and $Y$ is defined by

$$I(X,Y) = \mathcal{H}(P_{X,Y} \| P_X \otimes P_Y).$$

Now let $Z$ be a third random variable. The *conditional mutual information* of $X$ and $Y$ given $Z$ is defined by

$$I(X,Y|Z) = \mathbb{E}[\mathcal{H}(P_{X,Y|Z} \| P_{X|Z} \otimes P_{Y|Z})],$$

provided the regular conditional probabilities exist.

If $\mathcal{A}$ is a sub-$\sigma$-algebra of $\mathcal{F}$, then we write $\mathrm{id}_\mathcal{A}$ for the measurable map $(\Omega, \mathcal{F}) \to (\Omega, \mathcal{A}), \omega \mapsto \omega$. For two sub-$\sigma$-algebras $\mathcal{A}$ and $\mathcal{D}$, we abbreviate

$$I(\mathcal{A}, \mathcal{D}) = I(\mathrm{id}_\mathcal{A}, \mathrm{id}_\mathcal{D}).$$

Since our probability space is standard, for any sub-$\sigma$-fields $\mathcal{A}, \mathcal{D}, \mathcal{E}$ of $\mathcal{F}$, there exists a regular conditional probability $P_{\mathrm{id}_\mathcal{A}, \mathrm{id}_\mathcal{D} | \mathrm{id}_\mathcal{E}}$ and we define

$$I(\mathcal{A}, \mathcal{D} | \mathcal{E}) := I(\mathrm{id}_\mathcal{A}, \mathrm{id}_\mathcal{D} | \mathrm{id}_\mathcal{E}).$$

The mutual information was introduced by Shannon as a measure of information. It plays an important role in information theory (see, e.g., [16]).

THEOREM 5.8.

$$\lim_{|\Delta| \to 0} \sum_i I(\mathcal{G}_{s_i}^0, \mathcal{F}_{s_{i+1}}^0 | \mathcal{F}_{s_i}^0) = \tfrac{1}{2} E \int_0^T \mu_u^2 \, d\langle M, M \rangle_u.$$

PROOF. Note that, for three random variables $X, Y$ and $Z$, we have

$$\frac{dP_{(X,Y)|Z}}{d(P_{X|Z} \otimes P_{Y|Z})} = \frac{dP_{X|(Y,Z)}}{dP_{X|Z}}.$$

This property implies that one has, for $0 \leq s < t \leq T$,

$$I(\mathcal{G}_s^0, \mathcal{F}_t^0 | \mathcal{F}_s^0) = \iint \log \frac{dP_{\mathrm{id}_{\mathcal{G}_s^0} | \mathrm{id}_{\mathcal{F}_t^0}}}{dP_{\mathrm{id}_{\mathcal{G}_s^0} | \mathrm{id}_{\mathcal{F}_s^0}}} dP(\omega') dP(\omega)$$

$$= \iint \log \frac{P_t(\cdot, d\omega')}{P_s(\cdot, d\omega')} \bigg|_{\mathcal{G}_s^0} dP(\omega') dP(\omega)$$

$$= H_{\mathcal{G}_s^0}(s, t).$$

Thus, the assertion is an immediate consequence of Theorem 5.5. □

This result motivates the following notion.

DEFINITION 5.9. The *information difference* of $(\mathcal{G}_r^0)$ relative to $(\mathcal{F}_r^0)$ up to time $T$ is defined as

$$A(\mathcal{G}^0, \mathcal{F}^0) = \lim_{|\Delta| \to 0} \sum_i I(\mathcal{G}_{s_i}^0, \mathcal{F}_{s_{i+1}}^0 | \mathcal{F}_{s_i}^0).$$

REMARK 5.10. Note that we did not use $M$ in our definition of the information difference of $(\mathcal{G}_r^0)$ relative to $(\mathcal{F}_r^0)$. However, by Theorem 5.8, the information difference may be represented in terms of any local martingale satisfying the (PRP).

Theorem 5.8 can be reformulated in the following way.

THEOREM 5.11. *The additional utility of an agent with information* $(\mathcal{G}_t)$ *is equal to the information difference of* $(\mathcal{G}_r^0)$ *relative to* $(\mathcal{F}_r^0)$, *that is,*

$$u_{\mathcal{G}}(x) - u_{\mathcal{F}}(x) = A(\mathcal{G}^0, \mathcal{F}^0).$$

If $(\mathcal{G}_t)$ is initially enlarged by some random variable $G$, then the information difference of $(\mathcal{G}_r^0)$ relative to $(\mathcal{F}_r^0)$ coincides with the Shannon information between $G$ and $(\mathcal{F}_T^0)$.

LEMMA 5.12. *Let* $\mathcal{G}_t^0 = \mathcal{F}_t^0 \vee \sigma(G)$, *where* $G$ *is a random variable with values in some Polish space. Then*

$$A(\mathcal{G}^0, \mathcal{F}^0) = I(G, \mathcal{F}_T^0 | \mathcal{F}_0^0).$$

PROOF. Let $0 \leq s \leq t$. By standard arguments, we have $I(\mathcal{G}_s^0, \mathcal{F}_t^0 | \mathcal{F}_s^0) = I(G, \mathcal{F}_t^0 | \mathcal{F}_s^0)$ and

$$I(G, \mathcal{F}_t^0 | \mathcal{F}_0^0) = I(G, (\mathcal{F}_t^0, \mathcal{F}_s^0) | \mathcal{F}_0^0)$$

$$= I(G, \mathcal{F}_t^0 | \mathcal{F}_s^0) + I(G, \mathcal{F}_s^0 | \mathcal{F}_0^0)$$

(see, e.g., [16], Theorem 1.6.3). By iteration, we obtain, for all partitions $\Delta$,
$$\sum_i I(\mathcal{G}^0_{s_i}, \mathcal{F}^0_{s_{i+1}} | \mathcal{F}^0_{s_i}) = I(G, \mathcal{F}^0_T | \mathcal{F}^0_0)$$
and, hence, the result. $\square$

THEOREM 5.13. *Let $\mathcal{G}^0_t = \mathcal{F}^0_t \vee \sigma(G)$, where $G$ is a random variable with values in some Polish space. Then the additional logarithmic utility of an agent with information $(\mathcal{G}_t)$ is equal to the Shannon information between $G$ and $(\mathcal{F}^0_T)$ conditioned on $\mathcal{F}_0$, that is,*
$$u_\mathcal{G}(x) - u_\mathcal{F}(x) = I(\mathcal{F}^0_T, G | \mathcal{F}^0_0).$$
*In particular, if $\mathcal{F}^0_0$ is trivial, then the additional utility is equal to $I(\mathcal{F}^0_T, G)$.*

PROOF. This follows from Lemma 5.12 and Theorem 5.8. $\square$

REMARK 5.14. If $\mathcal{G}^0_t = \mathcal{F}^0_t \vee \sigma(G)$ and $G$ is $\mathcal{F}^0_T$-measurable, then the mutual information $I(\mathcal{F}^0_T, G | \mathcal{F}^0_0)$ is equal to the conditional absolute entropy of $G$ (see also [3]).

EXAMPLE 5.15. Let $(\Omega, \mathcal{F}, P)$ be the one-dimensional canonical Wiener space equipped with the Wiener process $(W_t)_{0 \leq t \leq 1}$. More precisely, $\Omega = \mathcal{C}([0,1], \mathbb{R})$ is the set of continuous functions on $[0,1]$ starting in 0, $\mathcal{F}$ the $\sigma$-algebra of Borel sets with respect to uniform convergence, $P$ the Wiener measure and $W$ the coordinate process. $(\mathcal{F}_t)_{0 \leq t \leq 1}$ is obtained by completing the natural filtration $(\mathcal{F}^0_t)_{0 \leq t \leq 1}$. Suppose the price process $S$ is of the form
$$S_t = \exp(W_t + bt), \qquad 0 \leq t \leq 1,$$
with $b \in \mathbb{R}$. We want to calculate the additional utility of an insider knowing whether the price exceeds a certain level or not. More precisely, we suppose the insider to know the value of
$$G = \mathbb{1}_{(c,\infty)}(S^*_1),$$
where $c > 0$ and $S^*_1 = \max_{0 \leq t \leq 1} S_t$. By Remark 5.14, the additional utility is equal to the entropy
$$H(G) = p \log p + (1-p) \log(1-p),$$
where
$$p = P(S^*_1 > c).$$
This may be calculated via Girsanov's theorem. Namely, we have
$$P(S^*_1 > c) = P\bigg(\forall t \in [0,1]: \max_{t \in [0,1]} W_t + bt > \log c\bigg)$$
$$= \int_0^1 \exp\bigg(b \log c - \frac{b^2}{2}s\bigg) \frac{|\log c|}{\sqrt{2\pi s^3}} \exp\bigg(-\frac{|\log c|^2}{2s}\bigg) ds.$$

**6. Mutual information estimates.** In this final section we apply some results from information theory to derive estimates for the information of a better informed agent. This yields a priori estimates for the agent's additional expected logarithmic utility in the light of the preceding section. Among other facts, the differential entropy maximizing property of Gaussian laws will play a role. We adopt the notation of [16].

Before we provide the information estimates, we summarize some basic facts of the mutual information (see [16], Theorem 1.6.3). For random variables $X, Y, Z$ in some Borel spaces, the following properties hold:

(I.1) $I(X, Y|Z) \geq 0$ and, $I(X, Y|Z) = 0$ if and only if $X$ and $Y$ are independent given $Z$.
(I.2) $I(X, (Y, Z)) = I(X, Z) + I(X, Y|Z)$.
(I.3) If $X$ is a continuous random variable with finite differential entropy, then
$$I(X, Y) = h(X) - h(X|Y).$$

For some fixed integer $d \in \mathbb{N}$, let $X$ be a $\mathcal{F}_T^0$-measurable $\mathbb{R}^d$-valued random variable. Moreover, denote by $Y$ a $d$-dimensional r.v. that is independent of the $\sigma$-field $\mathcal{F}_T^0$. We consider the enlarged filtration $\mathcal{G}_t^0 = \mathcal{F}_t^0 \vee \sigma(G)$, where $G := X + Y$.

LEMMA 6.1. *Suppose that the law of $Y$ is absolutely continuous with respect to Lebesgue measure and has finite differential entropy*
$$h(Y) = -\int \frac{dP_Y}{d\lambda^d}(y) \log \frac{dP_Y}{d\lambda^d}(y)\, dy.$$

*Then*

(12) $$I(G, \mathcal{F}_T^0) = h(X + Y) - h(Y).$$

PROOF. Due to property (I.2), we have
$$I(G, \mathcal{F}_T^0) = I(X + Y, X) + I(X + Y, \mathcal{F}_T^0 | X).$$
Given $X$, the r.v.'s $X + Y$ and $\mathrm{id}_{\mathcal{F}_T^0}$ are independent. Therefore, (I.1) and (I.3) lead to
$$I(G, \mathcal{F}_T^0) = I(X + Y, X) = h(X + Y) - h(X + Y|X) = h(X + Y) - h(Y).\ \square$$

Now assume the perturbation $Y$ to be a $\mathbb{R}^d$-valued centered Gaussian r.v. that is independent of $\mathcal{F}_T^0$.

LEMMA 6.2. *Suppose that $X \in L^2(P)$ and let $C_X$ and $C_Y$ denote the covariance matrices of $X$ and $Y$, respectively. Then*

$$I(G, \mathcal{F}_T^0) \leq \frac{1}{2} \log \frac{\det(C_X + C_Y)}{\det(C_Y)}. \tag{13}$$

*Moreover, equality holds in* (13) *if $X$ is Gaussian.*

PROOF. The distribution of $Y$ is continuous with respect to Lebesgue measure and has finite entropy. Therefore,

$$I(G, \mathcal{F}_T^0) = h(X+Y) - h(Y).$$

Let $C_X$ and $C_Y$ denote the covariance matrices of $X$ and $Y$, respectively. Due to the independence of $X$ and $Y$, the random variable $X+Y$ has the covariance matrix $C_{X+Y} = C_X + C_Y$. Next recall that the normal distribution maximizes the differential entropy under a covariance constraint, that is, $h(X+Y) \leq h(Z)$, where $Z$ is a centered Gaussian r.v. with covariance matrix $C_{X+Y}$. Therefore,

$$I(X, X+Y) \leq h(Z) - h(Y).$$

Using the formula for the differential entropy of Gaussian measures (Theorem 1.8.1, [16]), we obtain

$$h(Z) - h(Y) = \frac{1}{2} \log((2\pi e)^d \det(C_{X+Y})) - \frac{1}{2} \log((2\pi e)^d \det(C_Y))$$

$$= \frac{1}{2} \log \frac{\det(C_{X+Y})}{\det(C_Y)}.$$

If $X$ is Gaussian, then $h(X+Y) = h(Z)$ and, hence, the second statement of the lemma follows. □

COROLLARY 6.3. *Assume that additionally to the assumptions of the above lemma, the equation $Y = \kappa N$ is valid, where $N$ is a d-dimensional standard normal r.v. and $\kappa > 0$. Then*

$$I(G, \mathcal{F}_T^0) \leq \frac{1}{2} \sum_{j=1}^d \log \frac{\lambda_j + \kappa}{\kappa},$$

*where $\lambda_j$ ($j = 1, \ldots, n$) denote the eigenvalues of $C_X$.*

PROOF. The proof follows easily by computing the determinants in Lemma 6.2. □

The proof of Lemma 6.2 is based on the fact that Gaussian distributions maximize the differential entropy under a constraint on the covariance structure. Let us recall the construction of entropy maximizing measures under a linear constraint.

LEMMA 6.4. *Let $E \subset \mathbb{R}^d$ be a measurable set, $c > 0$ and $g : E \to [0, \infty)$ a measurable map. Assume that there exist constants $Z, t \geq 0$, such that the measure $\nu$ defined by*

$$\frac{d\nu}{d\lambda^d}(x) = \frac{1}{Z} e^{-tg(x)}$$

*is a probability measure satisfying $E^\nu[g] = c$. Then $\nu$ is the unique probability measure maximizing the differential entropy among all continuous probability measures $\mu$ on $E$ satisfying $E^\mu[g] = c$.*

The entropy maximization problem is equivalent to minimizing the relative entropy $\mathcal{H}(\cdot \| \lambda^d)$. Hence, the problem can be treated under more general constraints by using results of [9], Theorem 3.1.

PROOF. Let $\mu$ be a continuous probability measure on $E$ with $E^\mu[g] = c$. Then

$$\mathcal{H}(\mu \| \nu) = E^\mu \log \frac{d\mu}{d\nu} = E^\mu \log \frac{d\mu}{d\lambda^d} + E^\mu \log \frac{d\lambda^d}{d\nu}$$
$$= -h(\mu) + \log Z + t E^\mu g = -h(\mu) + \log Z + t E^\nu g$$
$$= -h(\mu) - E^\nu \log \frac{e^{-tg}}{Z} = -h(\mu) + h(\nu).$$

Since $\mathcal{H}(\mu \| \nu) \geq 0$ and $H(\mu \| \nu) = 0$ iff $\mu = \nu$, $\nu$ is the unique maximizer of the differential entropy. □

REMARK 6.5. The above lemma can be used to derive similar results as obtained in Lemma 6.2. For instance, for $E := \mathbb{R}$ and $g(x) := |x|$, one obtains that the two-sided exponential distribution maximizes the differential entropy under the constraint $E^\mu g = c$ ($c > 0$). In particular, the measure $\nu$ with $\frac{d\nu}{d\lambda}(x) = (2c)^{-1} e^{-|x|/c}$ satisfies

$$E^\nu[g] = c \quad \text{and} \quad h(\nu) = 1 + \log(2c).$$

Now let $X$ be a real-valued r.v. in $L^1(P)$. Moreover, let $\kappa_1 := E[|X - EX|]$ and $Y$ be a two-sided exponential distribution with $E|Y| =: \kappa_2$. Then due to Lemma 6.1,

$$I(G, \mathcal{F}_T^0) \leq \log\left(\frac{\kappa_1 + \kappa_2}{\kappa_2}\right).$$

EXAMPLE 6.6. We consider the classical stock market model with one asset. Let $(\mathcal{F}_t^0)_{t \in [0,T]}$ be a Brownian filtration generated by the Brownian motion $(B_t)_{t \in [0,T]}$ and denote by $(\mathcal{F}_t)$ its completion. The stock price is modeled by the process

$$S_t = S_0 \exp\{B_t + bt\},$$

where $S_0 > 0$ is the deterministic stock price at time 0 and $b \in \mathbb{R}$. For some fixed times $t_1, \ldots, t_d \in (0, T]$ ($d \in \mathbb{N}$), let $X := (B_{t_i})_{i=1,\ldots,d}$. We suppose that the insider bases his investment on the filtration $\mathcal{G}_t = \bigcap_{s>t} \mathcal{F}_s \vee \sigma(G)$, where $G = X + \kappa N$ and $N$ is a standard normal r.v. in $\mathbb{R}^d$ that is independent of $\mathcal{F}_T$. Due to Lemma 6.2, the additional utility of the insider is related to the eigenvalues of the matrix

$$\begin{pmatrix} t_1 & t_1 & \ldots & t_1 \\ t_1 & t_2 & \ldots & t_2 \\ \vdots & \vdots & & \vdots \\ t_1 & t_2 & \ldots & t_d \end{pmatrix}.$$

Let us finish the section with an example for a general enlargement.

EXAMPLE 6.7. We reconsider the classical stock market model of Example 6.6 with $T := 1$. The knowledge of the insider at time $t$ is modeled by $\mathcal{G}_t = \bigcap_{r>t} \mathcal{F}_r \vee \sigma((G_s)_{s \in [0,r]})$, where $G_t := B_1 + \widetilde{B}_{g(1-t)}$, $(\widetilde{B}_t)$ is a Brownian motion independent of $(B_t)$ and $g:[0,1] \to [0,\infty)$ is a decreasing function. We are therefore in a setting similar to Example 2.13. We now calculate the utility increment from the perspective of the notion of information difference of filtrations. Let $\pi$ be as in Section 4. For $0 \leq s \leq t \leq 1$, we have

$$\pi([0,s) \times (s,t]) = I((G_u)_{u \in [0,s]}, \mathcal{F}_t^0 | \mathcal{F}_s^0)$$
$$= I(G_s, \mathcal{F}_t^0 | \mathcal{F}_s^0) = I(G_s, B_t | \mathcal{F}_s^0) + I(G_s, \mathcal{F}_t^0 | \mathcal{F}_s^0, B_t)$$
$$= I(B_1 + \widetilde{B}_{g(1-s)}, B_t - B_s | \mathcal{F}_s^0)$$
$$= I(B_1 - B_s + \widetilde{B}_{g(1-s)}, B_t - B_s).$$

Using the formula for the differential entropy for Gaussian measures, we obtain

$$\pi([0,s) \times (s,t]) = h(B_1 - B_s + \widetilde{B}_{g(1-s)}) - h(B_1 - B_t + \widetilde{B}_{g(1-s)})$$
$$= \frac{1}{2} \log(2\pi e(1 - s + g(1-s)))$$
$$\quad - \frac{1}{2} \log(2\pi e(1 - t + g(1-s)))$$
$$= \frac{1}{2} \log \frac{1 - s + g(1-s)}{1 - t + g(1-s)}.$$

Alternatively, one can express $\pi([0,s) \times (s,t])$ as

$$\pi([0,s) \times (s,t]) = \frac{1}{2} \int_s^t \frac{1}{1 - u + g(1-s)} \, du.$$

For a partition $\Delta: 0 = t_0 \leq \cdots \leq t_m = 1$ ($m \in \mathbb{N}$), we consider $D^\Delta$ as in Section 4. One has

$$\pi(D^\Delta) = \sum_{i=1}^n \pi([t_{i-1}, t_i) \times (t_i, t_{i+1}])$$

$$= \frac{1}{2} \int_0^1 \frac{1}{1 - u + g(1 - \max\{t_i : t_i \leq u\})} \, du.$$

Next, choose a sequence of refining partitions $(\Delta_n)$ such that their mesh tends to 0. Then the term in the latter integral is monotonically increasing in $n$ and convergent. Hence, one obtains

$$\lim_{n \to \infty} \pi(D^{\Delta_n}) = \frac{1}{2} \int_0^1 \frac{1}{1 - u + g(1 - u)} \, du.$$

On the other hand,

$$\lim_{n \to \infty} \pi(D^{\Delta_n}) = \pi(D) = u_\mathcal{G}(x) - u_\mathcal{F}(x).$$

Consequently, the insider has finite utility if and only if $\int_0^1 \frac{1}{1-u+g(1-u)} \, du < \infty$ is finite. Now suppose $g(y) = Cy^p$ for some $C > 0$ and $p > 0$. It is straightforward to show that the integral and, hence, the additional utility, is finite if and only if $p \in (0, 1)$. This equivalence follows also from results in [8], where the authors compute explicitly the information drift.

S. Ankirchner
P. Imkeller
Institut für Mathematik
Humboldt-Universität zu Berlin
Unter den Linden 6
10099 Berlin
Germany
E-mail: imkeller@mathematik.hu-berlin.de

S. Dereich
Fachbereich Mathematik
Technische Universität Berlin
Strasse des 17. Juni 136
10623 Berlin
Germany
E-mail: dereich@math.tu-berlin.de